\documentclass{article}
\usepackage{amsfonts}
\usepackage{amsmath}
\usepackage{amssymb}
\usepackage{graphicx}

\newtheorem{theorem}{Theorem}

\newtheorem{lemma}[theorem]{Lemma}

\newtheorem{proposition}[theorem]{Proposition}
\newtheorem{remark}[theorem]{Remark}

\newenvironment{proof}[1][Proof]{\noindent\textbf{#1.} }{\ \rule{0.5em}{0.5em}}
\input{tcilatex}

\begin{document}

\title{Confidence Sets Based on Penalized Maximum Likelihood Estimators in
Gaussian Regression\thanks{%
Earlier versions of this paper were circulated under the title "Confidence
Sets Based on Penalized Maximum Likelihood Estimators".}}
\author{Benedikt M. P\"{o}tscher\thanks{%
Department of Statistics, University of Vienna, Universit\"{a}tsstrasse 5,
A-1010 Vienna. Phone: +431 427738640. E-mail: benedikt.poetscher@univie.ac.at%
} and Ulrike Schneider\thanks{%
Institute for Mathematical Stochastics, Georg-August-University G\"{o}%
ttingen, Goldschmidtstra\ss e 7, D-37077 G\"{o}ttingen. Phone: +49
55139172107. E-mail: ulrike.schneider@math.uni-goettingen.de} \\
Department of Statistics, University of Vienna\\
and\\
Institute for Mathematical Stochastics, University of G\"{o}ttingen}
\date{Preliminary version: February 2008\\
First version: June 2008\\
First revision: May 2009\\
Second revision: January 2010}
\maketitle

\begin{abstract}
Confidence intervals based on penalized maximum likelihood estimators such
as the LASSO, adaptive LASSO, and hard-thresholding are analyzed. In the
known-variance case, the finite-sample coverage properties of such intervals
are determined and it is shown that symmetric intervals are the shortest.
The length of the shortest intervals based on the hard-thresholding
estimator is larger than the length of the shortest interval based on the
adaptive LASSO, which is larger than the length of the shortest interval
based on the LASSO, which in turn is larger than the standard interval based
on the maximum likelihood estimator. In the case where the penalized
estimators are tuned to possess the `sparsity property', the intervals based
on these estimators are larger than the standard interval by an order of
magnitude. Furthermore, a simple asymptotic confidence interval construction
in the `sparse' case, that also applies to the smoothly clipped absolute
deviation estimator, is discussed. The results for the known-variance case
are shown to carry over to the unknown-variance case in an appropriate
asymptotic sense.

\emph{MSC Subject Classifications: }Primary 62F25; secondary 62C25,

62J07.

\emph{Keywords}: penalized maximum likelihood, penalized least squares,
Lasso, adaptive Lasso, hard-thresholding, soft-thresholding, confidence set,
coverage probability, sparsity, model selection.
\end{abstract}

\section{Introduction}

Recent years have seen an increased interest in penalized maximum likelihood
(least squares) estimators. Prominent examples of such estimators are the
LASSO estimator (Tibshirani (1996)) and its variants like the adaptive LASSO
(Zou (2006)), the Bridge estimators (Frank and Friedman (1993)), or the
smoothly clipped absolute deviation (SCAD) estimator (Fan and Li (2001)). In
linear regression models with orthogonal regressors, the hard- and
soft-thresholding estimators can also be reformulated as penalized least
squares estimators, with the soft-thresholding estimator then coinciding
with the LASSO estimator.

The asymptotic distributional properties of penalized maximum likelihood
(least squares) estimators have been studied in the literature, mostly in
the context of a finite-dimensional linear regression model; see Knight and
Fu (2000), Fan and Li (2001), and Zou (2006). Knight and Fu (2000) study the
asymptotic distribution of Bridge estimators and, in particular, of the
LASSO estimator. Their analysis concentrates on the case where the
estimators are tuned in such a way as to perform conservative model
selection, and their asymptotic framework allows for dependence of
parameters on sample size. In contrast, Fan and Li (2001) for the SCAD
estimator and Zou (2006) for the adaptive LASSO estimator concentrate on the
case where the estimators are tuned to possess the `sparsity' property. They
show that, with such tuning, these estimators possess what has come to be
known as the `oracle property'. However, their results are based on a
fixed-parameter asymptotic framework only. P\"{o}tscher and Leeb (2009) and P%
\"{o}tscher and Schneider (2009) study the finite-sample distribution of the
hard-thresholding, the soft-thresholding (LASSO), the SCAD, and the adaptive
LASSO estimator under normal errors; they also obtain the asymptotic
distributions of these estimators in a general `moving parameter' asymptotic
framework. The results obtained in these two papers clearly show that the
distributions of the estimators studied are often highly non-normal and that
the so-called `oracle property' typically paints a misleading picture of the
actual performance of the estimator. [In the wake of Fan and Li (2001) a
considerable literature has sprung up establishing the so-called `oracle
property' for a variety of estimators. All these results are fixed-parameter
asymptotic results only and can be very misleading. See Leeb and P\"{o}%
tscher (2008) and P\"{o}tscher (2009) for more discussion.]

A natural question now is what all these distributional results mean for
confidence intervals that are based on penalized maximum likelihood (least
squares) estimators. This is the question we address in the present paper in
the context of a normal linear regression model with orthogonal regressors.
In the known-variance case we obtain formulae for the finite-sample infimal
coverage probabilities of fixed-width confidence intervals based on the
following estimators: hard-thresholding, LASSO (soft-thresholding), and
adaptive LASSO. We show that among those intervals the symmetric ones are
the shortest, and we show that hard-thresholding leads to longer intervals
than the adaptive LASSO, which in turn leads to longer intervals than the
LASSO. All these intervals are longer than the standard confidence interval
based on the maximum likelihood estimator, which is in line with Joshi
(1969). In case the estimators are tuned to possess the `sparsity' property,
explicit asymptotic formulae for the length of the confidence intervals are
furthermore obtained, showing that in this case the intervals based on the
penalized maximum likelihood estimators are larger by an order of magnitude
than the standard maximum likelihood based interval. This refines, for the
particular estimators considered, a general result for confidence sets based
on `sparse' estimators (P\"{o}tscher (2009)). Additionally, in the
`sparsely' tuned case a simple asymptotic construction of confidence
intervals is provided that also applies to other penalized maximum
likelihood estimators such as the SCAD estimator. Furthermore, we show how
the results for the known-variance case carry over to the unknown-variance
case in an asymptotic sense.

The plan of the paper is as follows: After introducing the model and
estimators in Section 2, the known-variance case is treated in Section 3
whereas the unknown-variance case is dealt with in Section 4. All proofs as
well as some technical lemmata are relegated to the Appendix.

\section{The Model and Estimators}

For a normal linear regression model with orthogonal regressors,
distributional properties of penalized maximum likelihood (least squares)
estimators with a separable penalty can be reduced to the case of a Gaussian
location problem; for details see, e.g., P\"{o}tscher and Schneider (2009).
Since we are only interested in confidence sets for individual components of
the parameter vector in the regression that are based on such estimators, we
shall hence suppose that the data $y_{1},\ldots ,y_{n}$ are independent
identically distributed as $N(\theta ,\sigma ^{2})$, $\theta \in \mathbb{R}$%
, $0<\sigma <\infty $. [This entails no loss of generality in the
known-variance case. In the unknown-variance case an explicit treatment of
the orthogonal linear model would differ from the analysis in the present
paper only in that the estimator $\hat{\sigma}^{2}$ defined below would be
replaced by the usual residual variance estimator from the least-squares
regression; this would have no substantial effect on the results.] We shall
be concerned with confidence sets for $\theta $ based on penalized maximum
likelihood estimators such as the hard-thresholding estimator, the LASSO
(reducing to soft-thresholding in this setting), and the adaptive LASSO
estimator. The hard-thresholding estimator $\tilde{\theta}_{H}$ is given by%
\begin{equation*}
\tilde{\theta}_{H}:=\tilde{\theta}_{H}(\eta _{n})=\bar{y}\boldsymbol{1}%
(\left\vert \bar{y}\right\vert >\hat{\sigma}\eta _{n})
\end{equation*}%
where the threshold $\eta _{n}$ is a positive real number, $\bar{y}$ denotes
the maximum likelihood estimator, i.e., the arithmetic mean of the data, and 
$\hat{\sigma}^{2}=(n-1)^{-1}\sum_{i=1}^{n}(y_{i}-\bar{y})^{2}$. Also define
the infeasible estimator%
\begin{equation*}
\hat{\theta}_{H}:=\hat{\theta}_{H}(\eta _{n})=\bar{y}\boldsymbol{1}%
(\left\vert \bar{y}\right\vert >\sigma \eta _{n})
\end{equation*}%
which uses the value of $\sigma $. The LASSO (or soft-thresholding)
estimator $\tilde{\theta}_{S}$ is given by 
\begin{equation*}
\tilde{\theta}_{S}:=\tilde{\theta}_{S}(\eta _{n})=\limfunc{sign}(\bar{y}%
)(\left\vert \bar{y}\right\vert -\hat{\sigma}\eta _{n})_{+}
\end{equation*}%
and its infeasible version by%
\begin{equation*}
\hat{\theta}_{S}:=\hat{\theta}_{S}(\eta _{n})=\limfunc{sign}(\bar{y}%
)(\left\vert \bar{y}\right\vert -\sigma \eta _{n})_{+}.
\end{equation*}%
Here $\limfunc{sign}(x)$ is defined as $-1$, $0$, and $1$ in case $x<0$, $%
x=0 $, and $x>0$, respectively, and $z_{+}$ is shorthand for $\max \{z,0\}$.
The adaptive LASSO estimator $\tilde{\theta}_{A}$ in this simple model is
given by

\begin{equation*}
\tilde{\theta}_{A}:=\tilde{\theta}_{A}(\eta _{n})=\bar{y}(1-\hat{\sigma}%
^{2}\eta _{n}^{2}/\bar{y}^{2})_{+}=\left\{ 
\begin{array}{cl}
0 & \text{if }\;|\bar{y}|\leq \hat{\sigma}\eta _{n} \\ 
\bar{y}-\hat{\sigma}^{2}\eta _{n}^{2}/\bar{y} & \text{if}\;\;|\bar{y}|>\hat{%
\sigma}\eta _{n},%
\end{array}%
\right.
\end{equation*}%
and its infeasible counterpart by%
\begin{equation*}
\hat{\theta}_{A}:=\hat{\theta}_{A}(\eta _{n})=\bar{y}(1-\sigma ^{2}\eta
_{n}^{2}/\bar{y}^{2})_{+}=\left\{ 
\begin{array}{cl}
0 & \text{if }\;|\bar{y}|\leq \sigma \eta _{n} \\ 
\bar{y}-\sigma ^{2}\eta _{n}^{2}/\bar{y} & \text{if}\;\;|\bar{y}|>\sigma
\eta _{n}.%
\end{array}%
\right.
\end{equation*}%
It coincides with the nonnegative Garotte in this simple model. For the
feasible estimators we always need to assume $n\geq 2$, whereas for the
infeasible estimators also $n=1$ is admissible.

Note that $\eta _{n}$ plays the r\^{o}le of a tuning parameter and it is
most natural to let the estimators depend on the tuning parameter only via $%
\sigma \eta _{n}$ and $\hat{\sigma}\eta _{n}$, respectively, in order to
take account of the scale of the data. This makes the estimators mentioned
above scale equivariant. We shall often suppress dependence of the
estimators on $\eta _{n}$ in the notation. In the following let $P_{n,\theta
,\sigma }$ denote the distribution of the sample when $\theta $ and $\sigma $
are the true parameters. Furthermore, let $\Phi $ denote the standard normal
cumulative distribution function.

We also note the following obvious fact: Since hard- and soft-thresholding
operate in a coordinatewise fashion, the results given below also apply
mutatis mutandis to linear regressions with non-orthogonal regressors. Of
course, the soft-thresholding estimator then no longer coincides with the
LASSO estimator. We refrain from spelling out details.

\section{Confidence Intervals: Known-Variance Case\label{finite_sample}}

In this section we consider the case where the variance $\sigma ^{2}$ is
known, $n\geq 1$ holds, and we are interested in the finite-sample coverage
properties of intervals of the form $[\hat{\theta}-\sigma a_{n},\hat{\theta}%
+\sigma b_{n}]$ where $a_{n}$ and $b_{n}$ are nonnegative real numbers and $%
\hat{\theta}$ stands for any one of the estimators $\hat{\theta}_{H}=\hat{%
\theta}_{H}(\eta _{n})$, $\hat{\theta}_{S}=\hat{\theta}_{S}(\eta _{n})$, or $%
\hat{\theta}_{A}=\hat{\theta}_{A}(\eta _{n})$. We shall also consider
one-sided intervals $(-\infty ,\hat{\theta}+\sigma c_{n}]$ and $[\hat{\theta}%
-\sigma c_{n},\infty )$ with $0\leq c_{n}<\infty $. Let $p_{n}(\theta
;\sigma ,\eta _{n},a_{n},b_{n})=P_{n,\theta ,\sigma }\left( \theta \in
\lbrack \hat{\theta}-\sigma a_{n},\hat{\theta}+\sigma b_{n}]\right) $ denote
the coverage probability. Due to the above-noted scale equivariance of the
estimator $\hat{\theta}$, it is obvious that 
\begin{equation*}
p_{n}(\theta ;\sigma ,\eta _{n},a_{n},b_{n})=p_{n}(\theta /\sigma ;1,\eta
_{n},a_{n},b_{n})
\end{equation*}%
holds, and the same is true for the one-sided intervals. In particular, it
follows that the infimal coverage probabilities $\inf_{\theta \in \mathbb{R}%
}p_{n}(\theta ;\sigma ,\eta _{n},a_{n},b_{n})$ do not depend on $\sigma $.
Therefore, we shall assume without loss of generality that $\sigma =1$ for
the remainder of this section and we shall write $P_{n,\theta }$ for $%
P_{n,\theta ,1}$.

\subsection{Infimal coverage probabilities in finite samples\label%
{inf_cov_prob}}

We begin with soft-thresholding. Let $C_{S,n}$ denote the interval $[\hat{%
\theta}_{S}-a_{n},\hat{\theta}_{S}+b_{n}]$. We first determine the infimum
of the coverage probability $p_{S,n}(\theta ):=p_{S,n}(\theta ;1,\eta
_{n},a_{n},b_{n})=P_{n,\theta }\left( \theta \in C_{S,n}\right) $ of this
interval.

\begin{proposition}
\label{lasso} For every $n\geq 1$, the infimal coverage probability of the
interval $C_{S,n}$ is given by%
\begin{equation}
\inf_{\theta \in \mathbb{R}}p_{S,n}(\theta )=\left\{ 
\begin{array}{cc}
\Phi (n^{1/2}(a_{n}-\eta _{n}))-\Phi (n^{1/2}(-b_{n}-\eta _{n})) & \text{if
\ }a_{n}\leq b_{n} \\ 
\Phi (n^{1/2}(b_{n}-\eta _{n}))-\Phi (n^{1/2}(-a_{n}-\eta _{n})) & \text{if
\ }a_{n}>b_{n}.%
\end{array}%
\right.  \label{infimal_soft_asym}
\end{equation}
\end{proposition}

As a point of interest we note that $p_{S,n}(\theta )$ is a piecewise
constant function with jumps at $\theta =-a_{n}$ and $\theta =b_{n}$.

Next we turn to hard-thresholding. Let $C_{H,n}$ denote the interval $[\hat{%
\theta}_{H}-a_{n},\hat{\theta}_{H}+b_{n}]$. The infimum of the coverage
probability $p_{H,n}(\theta ):=p_{H,n}(\theta ;1,\eta
_{n},a_{n},b_{n})=P_{n,\theta }\left( \theta \in C_{H,n}\right) $ of this
interval has been obtained in Proposition 3.1 in P\"{o}tscher (2009), which
we repeat for convenience.

\begin{proposition}
\label{hard} For every $n\geq 1$, the infimal coverage probability of the
interval $C_{H,n}$ is given by%
\begin{eqnarray}
&&\inf_{\theta \in \mathbb{R}}p_{H,n}(\theta )  \label{infimal_hard_asymm} \\
&=&\left\{ 
\begin{array}{ll}
\Phi (n^{1/2}(a_{n}-\eta _{n}))-\Phi (-n^{1/2}b_{n}) & \text{if \ \ }\eta
_{n}\leq a_{n}+b_{n}\text{ \ and \ }a_{n}\leq b_{n} \\ 
\Phi (n^{1/2}(b_{n}-\eta _{n}))-\Phi (-n^{1/2}a_{n}) & \text{if \ \ }\eta
_{n}\leq a_{n}+b_{n}\text{ \ and \ }a_{n}>b_{n} \\ 
0 & \text{if \ \ }\eta _{n}>a_{n}+b_{n}.%
\end{array}%
\right.  \notag
\end{eqnarray}
\end{proposition}

For later use we observe that the interval $C_{H,n}$ has positive infimal
coverage probability if and only if the length of the interval $a_{n}+b_{n}$
is larger than $\eta _{n}$. As a point of interest we also note that the
coverage probability $p_{H,n}(\theta )$ is discontinuous (with discontinuity
points at $\theta =-a_{n}$ and $\theta =b_{n}$). Furthermore, as discussed
in P\"{o}tscher (2009), the infimum in (\ref{infimal_hard_asymm}) is
attained if $\eta _{n}>a_{n}+b_{n}$, but not in case $\eta _{n}\leq
a_{n}+b_{n}$.

Finally, we consider the adaptive LASSO. Let $C_{A,n}$ denote the interval $[%
\hat{\theta}_{A}-a_{n},\hat{\theta}_{A}+b_{n}]$. The infimum of the coverage
probability $p_{A,n}(\theta ):=p_{A,n}(\theta ;1,\eta
_{n},a_{n},b_{n})=P_{n,\theta }\left( \theta \in C_{A,n}\right) $ of this
interval is given next.

\begin{proposition}
\label{adLASSOinf} For every $n\geq 1$, the infimal coverage probability of $%
C_{A,n}$ is given by%
\begin{equation*}
\inf_{\theta \in \mathbb{R}}p_{A,n}(\theta )=\Phi (n^{1/2}(a_{n}-\eta
_{n}))-\Phi \left( n^{1/2}\left( (a_{n}-b_{n})/2-\sqrt{%
((a_{n}+b_{n})/2)^{2}+\eta _{n}^{2}}\right) \right)
\end{equation*}%
if $a_{n}\leq b_{n}$, and by%
\begin{equation*}
\inf_{\theta \in \mathbb{R}}p_{A,n}(\theta )=\Phi (n^{1/2}(b_{n}-\eta
_{n}))-\Phi \left( n^{1/2}\left( (b_{n}-a_{n})/2-\sqrt{%
((a_{n}+b_{n})/2)^{2}+\eta _{n}^{2}}\right) \right)
\end{equation*}%
if $a_{n}>b_{n}$.
\end{proposition}

We note that $p_{A,n}$ is continuous except at $\theta =b_{n}$ and $\theta
=-a_{n}$ and that the infimum of $p_{A,n}$ is not attained which can be seen
from a simple refinement of the proof of Proposition \ref{adLASSOinf}.

\begin{remark}
\label{open}(i) If we consider the open interval $C_{S,n}^{o}=(\hat{\theta}%
_{S}-a_{n},\hat{\theta}_{S}+b_{n})$ the formula for the coverage probability
becomes 
\begin{eqnarray*}
P_{n,\theta }\left( \theta \in C_{S,n}^{o}\right) &=&[\Phi
(n^{1/2}(a_{n}-\eta _{n}))-\Phi (n^{1/2}(-b_{n}-\eta _{n}))]\boldsymbol{1}%
(\theta \leq -a_{n}) \\
&+&[\Phi (n^{1/2}(a_{n}+\eta _{n}))-\Phi (n^{1/2}(-b_{n}-\eta _{n}))]%
\boldsymbol{1}(-a_{n}<\theta <b_{n}) \\
&+&[\Phi (n^{1/2}(a_{n}+\eta _{n}))-\Phi (n^{1/2}(-b_{n}+\eta _{n}))]%
\boldsymbol{1}(b_{n}\leq \theta ).
\end{eqnarray*}%
As a consequence, the infimal coverage probability of $C_{S,n}^{o}$ is again
given by (\ref{infimal_soft_asym}). A fortiori, the half-open intervals $(%
\hat{\theta}_{n}-a_{n},\hat{\theta}_{n}+b_{n}]$ and $[\hat{\theta}_{n}-a_{n},%
\hat{\theta}_{n}+b_{n})$ then also have infimal coverage probability given
by (\ref{infimal_soft_asym}).

(ii) For the open interval $C_{H,n}^{o}=(\hat{\theta}_{H}-a_{n},\hat{\theta}%
_{H}+b_{n})$ the coverage probability satisfies 
\begin{multline*}
P_{n,\theta }\left( \theta \in C_{H,n}^{o}\right) =P_{n,\theta }\left(
\theta \in C_{H,n}\right)  \\
-[\boldsymbol{1}(\theta =b_{n})+\boldsymbol{1}(\theta =-a_{n})][\Phi
(n^{1/2}(-\theta +\eta _{n}))-\Phi (n^{1/2}(-\theta -\eta _{n}))].
\end{multline*}%
Inspection of the proof of Proposition 3.1 in P\"{o}tscher (2009) then shows
that $C_{H,n}^{o}$ has the same infimal coverage probability as $C_{H,n}$.
However, now the infimum is always a minimum. Furthermore, the half-open
intervals $(\hat{\theta}_{H}-a_{n},\hat{\theta}_{H}+b_{n}]$ and $[\hat{\theta%
}_{H}-a_{n},\hat{\theta}_{H}+b_{n})$ then a fortiori have infimal coverage
probability given by (\ref{infimal_hard_asymm}); for these intervals the
infimum is attained if $\eta _{n}>a_{n}+b_{n}$, but not necessarily if $\eta
_{n}\leq a_{n}+b_{n}$.

(iii) If $C_{A,n}^{o}$ denotes the open interval $(\hat{\theta}_{A}-a_{n},%
\hat{\theta}_{A}+b_{n})$, the formula for the coverage probability becomes%
\begin{multline*}
P_{n,\theta }\left( \theta \in C_{A,n}^{o}\right) = \\
\left\{ 
\begin{array}{ll}
\Phi \left( n^{1/2}\gamma ^{(-)}(\theta ,-a_{n})\right) -\Phi \left(
n^{1/2}\gamma ^{(-)}(\theta ,b_{n})\right)  & \text{if }\;\theta \leq -a_{n}
\\ 
\Phi \left( n^{1/2}\gamma ^{(+)}(\theta ,-a_{n})\right) -\Phi \left(
n^{1/2}\gamma ^{(-)}(\theta ,b_{n})\right)  & \text{if }\;-a_{n}<\theta
<b_{n} \\ 
\Phi \left( n^{1/2}\gamma ^{(+)}(\theta ,-a_{n})\right) -\Phi \left(
n^{1/2}\gamma ^{(+)}(\theta ,b_{n})\right)  & \text{if }\;\theta \geq b_{n},%
\end{array}%
\right. 
\end{multline*}%
where $\gamma ^{(-)}$ and $\gamma ^{(+)}$ are defined in (\ref{gamma1}) and (%
\ref{gamma2}) in the Appendix. Again the coverage probability is continuous
except at $\theta =b_{n}$ and $\theta =-a_{n}$ (and is continuous everywhere
in the trivial case $a_{n}=b_{n}=0$). It is now easy to see that the infimal
coverage probability of $C_{A,n}^{o}$ coincides with the infimal coverage
probability of the closed interval $C_{A,n}$, the infimum of the coverage
probability of $C_{A,n}^{o}$ now always being a minimum. Furthermore, the
half-open intervals $(\hat{\theta}_{A}-a_{n},\hat{\theta}_{A}+b_{n}]$ and $[%
\hat{\theta}_{A}-a_{n},\hat{\theta}_{A}+b_{n})$ a fortiori have the same
infimal coverage probability as $C_{A,n}$ and $C_{A,n}^{o}$.

(iv) The one-sided intervals $(-\infty ,\hat{\theta}_{S}+c_{n}]$, $(-\infty ,%
\hat{\theta}_{S}+c_{n})$, $[\hat{\theta}_{S}-c_{n},\infty )$, $(\hat{\theta}%
_{S}-c_{n},\infty )$, $(-\infty ,\hat{\theta}_{H}+c_{n}]$, $(-\infty ,\hat{%
\theta}_{H}+c_{n})$, $[\hat{\theta}_{H}-c_{n},\infty )$, $(\hat{\theta}%
_{H}-c_{n},\infty )$, $(-\infty ,\hat{\theta}_{A}+c_{n}]$, $(-\infty ,\hat{%
\theta}_{A}+c_{n})$, $(\hat{\theta}_{A}-c_{n},\infty )$, and $[\hat{\theta}%
_{A}-c_{n},\infty )$, with $c_{n}$ a nonnegative real number, have infimal
coverage probability $\Phi (n^{1/2}(c_{n}-\eta _{n}))$. This is easy to see
for soft-thresholding, follows from the reasoning in P\"{o}tscher (2009) for
hard-thresholding, and for the adaptive LASSO\ follows by similar, but
simpler, reasoning as in the proof of Proposition \ref{adLASSOinf}.
\end{remark}

\subsection{Symmetric intervals are shortest\label{symm}}

For the two-sided confidence sets considered above, we next show that given
a prescribed infimal coverage probability the symmetric intervals are
shortest. We then show that these shortest intervals are longer than the
standard interval based on the maximum likelihood estimator and quantify the
excess length of these intervals.

\begin{theorem}
\label{short}For every $n\geq 1$ and every $\delta $ satisfying $0<\delta <1$
we have:

(a) Among all intervals $C_{S,n}$ with infimal coverage probability not less
than $\delta $ there is a unique shortest interval $C_{S,n}^{\ast }=[\hat{%
\theta}_{S}-a_{n,S}^{\ast },\hat{\theta}_{S}+b_{n,S}^{\ast }]$ characterized
by $a_{n,S}^{\ast }=b_{n,S}^{\ast }$ with $a_{n,S}^{\ast }$ being the unique
solution of 
\begin{equation}
\Phi (n^{1/2}(a_{n}-\eta _{n}))-\Phi (n^{1/2}(-a_{n}-\eta _{n}))=\delta .
\label{short_a_S}
\end{equation}%
The interval $C_{S,n}^{\ast }$\ has infimal coverage probability equal to $%
\delta $ and $a_{n,S}^{\ast }$ is positive.

(b) Among all intervals $C_{H,n}$ with infimal coverage probability not less
than $\delta $ there is a unique shortest interval $C_{H,n}^{\ast }=[\hat{%
\theta}_{H}-a_{n,H}^{\ast },\hat{\theta}_{H}+b_{n,H}^{\ast }]$ characterized
by $a_{n,H}^{\ast }=b_{n,H}^{\ast }$ with $a_{n,H}^{\ast }$ being the unique
solution of 
\begin{equation}
\Phi (n^{1/2}(a_{n}-\eta _{n}))-\Phi (-n^{1/2}a_{n})=\delta .
\label{short_a_H}
\end{equation}%
The interval $C_{H,n}^{\ast }$ has infimal coverage probability equal to $%
\delta $ and $a_{n,H}^{\ast }$ satisfies $a_{n,H}^{\ast }>\eta _{n}/2$.

(c) Among all intervals $C_{A,n}$ with infimal coverage probability not less
than $\delta $ there is a unique shortest interval $C_{A,n}^{\ast }=[\hat{%
\theta}_{A}-a_{n,A}^{\ast },\hat{\theta}_{A}+b_{n,A}^{\ast }]$ characterized
by $a_{n,A}^{\ast }=b_{n,A}^{\ast }$ with $a_{n,A}^{\ast }$ being the unique
solution of%
\begin{equation}
\Phi (n^{1/2}(a_{n}-\eta _{n}))-\Phi \left( -n^{1/2}\sqrt{a_{n}^{2}+\eta
_{n}^{2}}\right) =\delta .  \label{short_a_A}
\end{equation}%
The interval $C_{A,n}^{\ast }$ has infimal coverage probability equal to $%
\delta $ and $a_{n,A}^{\ast }$ is positive.
\end{theorem}

In the statistically uninteresting case $\delta =0$ the interval with $%
a_{n}=b_{n}=0$ is the unique shortest interval in all three cases. However,
for the case of the hard-thresholding estimator also any interval with $%
a_{n}=b_{n}$ and $a_{n}\leq \eta _{n}/2$ has infimal coverage probability
equal to zero.

Given that the distributions of the estimation errors $\hat{\theta}%
_{S}-\theta $, $\hat{\theta}_{H}-\theta $, and $\hat{\theta}_{A}-\theta $
are not symmetric (see P\"{o}tscher and Leeb (2009), P\"{o}tscher and
Schneider (2009)), it may seem surprising at first glance that the shortest
confidence intervals are symmetric. Some intuition for this phenomenon can
be gained on the grounds that the distributions of the estimation errors
under $\theta =\tau $ and $\theta =-\tau $ are mirror-images of one another.

The above theorem shows that given a prespecified $\delta $ ($0<\delta <1$),
the shortest confidence set with infimal coverage probability equal to $%
\delta $ based on the soft-thresholding (LASSO) estimator is shorter than
the corresponding interval based on the adaptive LASSO estimator, which in
turn is shorter than the corresponding interval based on the
hard-thresholding estimator. All three intervals are longer than the
corresponding standard confidence interval based on the maximum likelihood
estimator. That is, 
\begin{equation*}
a_{n,H}^{\ast }>a_{n,A}^{\ast }>a_{n,S}^{\ast }>n^{-1/2}\Phi ^{-1}((1+\delta
)/2).
\end{equation*}%
Figure 1 below shows $n^{1/2}$ times the half-length of the shortest $\delta 
$-level confidence intervals based on hard-thresholding, adaptive LASSO,
soft-thresholding, and the maximum likelihood estimator, respectively, as a
function of $n^{1/2}\eta _{n}$ for various values of $\delta $. The graphs
illustrate that the intervals based on hard-thresholding, adaptive LASSO,
and soft-thresholding substantially exceed the length of the maximum
likelihood based interval except if $n^{1/2}\eta _{n}$ is very small. For
large values of $n^{1/2}\eta _{n}$ the graphs suggest a linear increase in
the length of the intervals based on the penalized estimators. This is
formally confirmed in Section \ref{asy_length} below.


\begin{figure}[p] 
\begin{center}
\includegraphics[width=\textwidth]{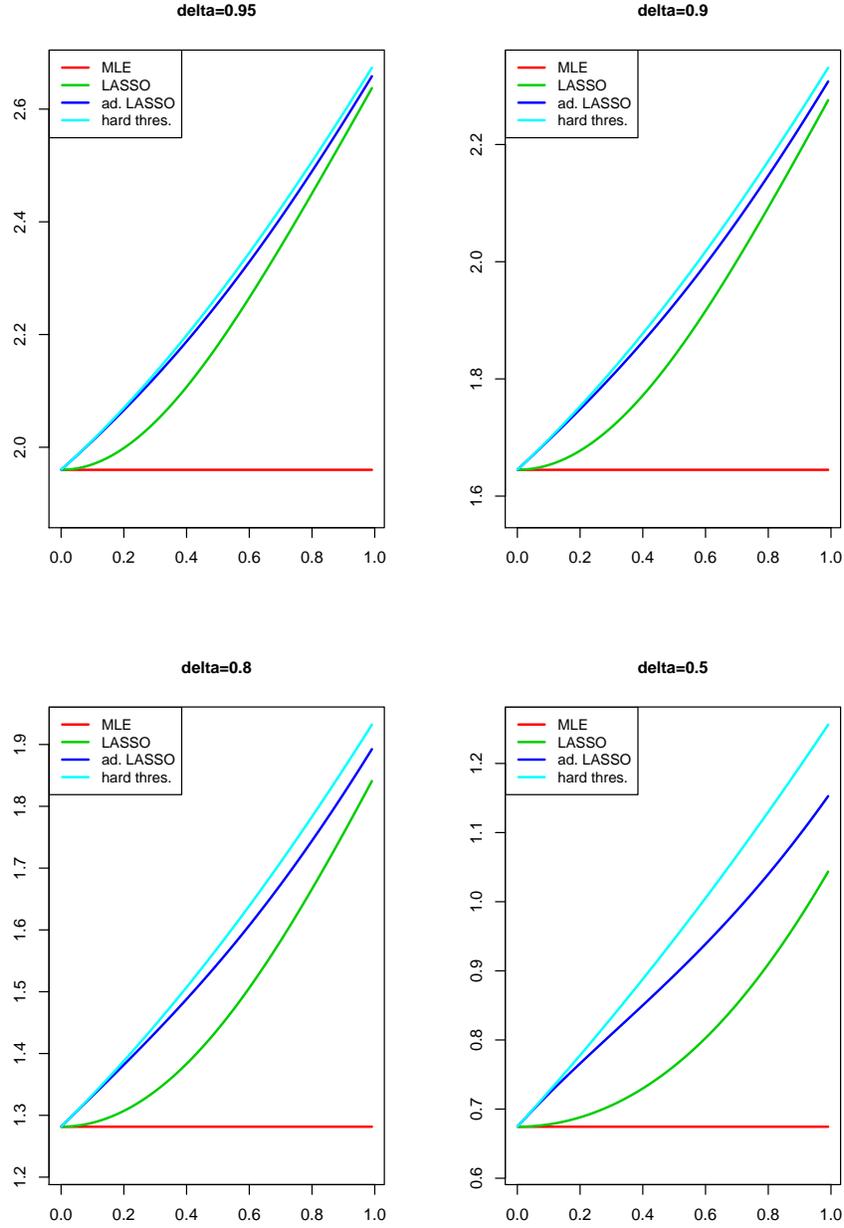}
\caption{$n^{1/2}a_{n,H}^{\ast }$, $n^{1/2}a_{n,A}^{\ast }$, 
$n^{1/2}a_{n,S}^{\ast }$ as a function of $n^{1/2}\eta _{n}$ for coverage 
probabilities $\protect\delta =0.5$, $0.8$, $0.9$, $0.95$. 
The horizontal line at height $\Phi^{-1}((1+\delta)/2)$ indicates $n^{1/2}$ 
times the half-length of the standard maximum likelihood based interval.}
\end{center}
\end{figure}

\newpage

\subsubsection{Asymptotic behavior of the length\label{asy_length}}

It is well-known that as $n\rightarrow \infty $ two different regimes for
the tuning parameter $\eta _{n}$ can be distinguished. In the first regime $%
\eta _{n}\rightarrow 0$ and $n^{1/2}\eta _{n}\rightarrow e$, $0<e<\infty $.
This choice of tuning parameter leads to estimators $\hat{\theta}_{S}$, $%
\hat{\theta}_{H}$, and $\hat{\theta}_{A}$ that perform conservative model
selection. In the second regime $\eta _{n}\rightarrow 0$ and $n^{1/2}\eta
_{n}\rightarrow \infty $, leading to estimators $\hat{\theta}_{S}$, $\hat{%
\theta}_{H}$, and $\hat{\theta}_{A}$ that perform consistent model selection
(also known as the `sparsity property');\ that is, with probability
approaching $1$, the estimators are exactly zero if the true value $\theta
=0 $, and they are different from zero if $\theta \neq 0$. See P\"{o}tscher
and Leeb (2009) and P\"{o}tscher and Schneider (2009) for a detailed
discussion. We now discuss the asymptotic behavior, under the two regimes,
of the half-length $a_{n,S}^{\ast }$, $a_{n,H}^{\ast }$, and $a_{n,A}^{\ast
} $ of the shortest intervals $C_{S,n}^{\ast }$, $C_{H,n}^{\ast }$, and $%
C_{A,n}^{\ast }$ with a fixed infimal coverage probability $\delta $, $%
0<\delta <1$.

If $\eta _{n}\rightarrow 0$ and $n^{1/2}\eta _{n}\rightarrow e$, $0<e<\infty 
$, then it follows immediately from Theorem \ref{short} that $%
n^{1/2}a_{n,S}^{\ast }$, $n^{1/2}a_{n,H}^{\ast }$, and $n^{1/2}a_{n,A}^{\ast
}$ converge to the unique solutions of 
\begin{equation}
\Phi (a-e)-\Phi (-a-e)=\delta ,  \label{asy_half-lenght_S_0}
\end{equation}%
\begin{equation}
\Phi (a-e)-\Phi (-a)=\delta ,  \label{asy_half-lenght_H_0}
\end{equation}%
and%
\begin{equation}
\Phi \left( \sqrt{a^{2}+e^{2}}\right) -\Phi (-a+e)=\delta ,
\label{asy_half-lenght_A_0}
\end{equation}%
respectively. [Actually, this is even true if $e=0$.] Hence, while $%
a_{n,H}^{\ast }$, $a_{n,A}^{\ast }$, and $a_{n,S}^{\ast }$ are larger than
the half-length $n^{-1/2}\Phi ^{-1}((1+\delta )/2)$ of the standard
interval, they are of the same order $n^{-1/2}$.

The situation is different, however, if $\eta _{n}\rightarrow 0$ but $%
n^{1/2}\eta _{n}\rightarrow \infty $. In this case Theorem \ref{short} shows
that 
\begin{equation*}
\Phi (n^{1/2}(a_{n,S}^{\ast }-\eta _{n}))\rightarrow \delta
\end{equation*}%
since $n^{1/2}(-a_{n,S}^{\ast }-\eta _{n})\leq -n^{1/2}\eta _{n}\rightarrow
-\infty $. In other words, 
\begin{equation}
a_{n,S}^{\ast }=\eta _{n}+n^{-1/2}\Phi ^{-1}(\delta )+o(n^{-1/2}).
\label{asy_half-length_S}
\end{equation}%
Similarly, noting that $n^{1/2}a_{n,H}^{\ast }>n^{1/2}\eta _{n}/2\rightarrow
\infty $, we get%
\begin{equation}
a_{n,H}^{\ast }=\eta _{n}+n^{-1/2}\Phi ^{-1}(\delta )+o(n^{-1/2});
\label{asy_half-length_H}
\end{equation}%
and since $n^{1/2}\sqrt{a_{n}^{2}+\eta _{n}^{2}}\geq n^{1/2}\eta
_{n}\rightarrow \infty $ we obtain%
\begin{equation}
a_{n,A}^{\ast }=\eta _{n}+n^{-1/2}\Phi ^{-1}(\delta )+o(n^{-1/2}).
\label{asy_half-length_A}
\end{equation}%
[Actually, the condition $\eta _{n}\rightarrow 0$ has not been used in the
derivation of (\ref{asy_half-length_S})-(\ref{asy_half-length_A}).] Hence,
the intervals $C_{S,n}^{\ast }$, $C_{H,n}^{\ast }$, and $C_{A,n}^{\ast }$
are asymptotically of the same length. They are also longer than the
standard interval by an order of magnitude: the ratio of each of $%
a_{n,S}^{\ast }$ ($a_{n,H}^{\ast }$, $a_{n,A}^{\ast }$, respectively) to the
half-length of the standard interval is $n^{1/2}\eta _{n}$, which diverges
to infinity. Hence, when the estimators $\hat{\theta}_{S}$, $\hat{\theta}%
_{H} $, and $\hat{\theta}_{A}$ are tuned to possess the `sparsity property',
the corresponding confidence sets become very large. For the particular
intervals considered here this is a refinement of a general result in P\"{o}%
tscher (2009) for confidence sets based on arbitrary estimators possessing
the `sparsity property'. [We note that the sparsely tuned hard-thresholding
estimator or the sparsely tuned adaptive LASSO (under an additional
condition on $\eta _{n}$) are known to possess the so-called `oracle
property'. In light of the `oracle property' it is sometimes argued in the
literature that valid confidence intervals based on these estimators with
length proportional to $n^{-1/2}$ can be obtained. However, in light of the
above discussion such intervals necessarily have infimal coverage
probability that converges to zero and thus are not valid. This once more
shows that \emph{fixed-parameter }asymptotic results like the `oracle'
property can be dangerously misleading.]

\subsection{A simple asymptotic confidence interval}

The results for the finite-sample confidence intervals given in Section \ref%
{inf_cov_prob} required a detailed case by case analysis based on the
finite-sample distribution of the estimator on which the interval is based.
If the estimators $\hat{\theta}_{S}$, $\hat{\theta}_{H}$, and $\hat{\theta}%
_{A}$ are tuned to possess the `sparsity property', i.e., if the tuning
parameter satisfies $\eta _{n}\rightarrow 0$ and $n^{1/2}\eta
_{n}\rightarrow \infty $, a simple asymptotic confidence interval
construction relying on asymptotic results obtained in P\"{o}tscher and Leeb
(2009) and P\"{o}tscher and Schneider (2009) is possible as shown below. An
advantage of this construction is that it easily extends to other estimators
like the smoothly clipped absolute deviation (SCAD) estimator when tuned to
possess the `sparsity property'.

As shown in P\"{o}tscher and Leeb (2009) and P\"{o}tscher and Schneider
(2009), the uniform rate of consistency of the `sparsely' tuned estimators $%
\hat{\theta}_{S}$, $\hat{\theta}_{H}$, and $\hat{\theta}_{A}$ is not $%
n^{1/2} $, but only $\eta _{n}^{-1}$; furthermore, the limiting
distributions of these estimators under the appropriate $\eta _{n}^{-1}$%
-scaling and under a moving-parameter asymptotic framework are always
concentrated on the interval $[-1,1]$. These facts can be used to obtain the
following result.

\begin{proposition}
\label{asy}Suppose $\eta _{n}\rightarrow 0$ and $n^{1/2}\eta _{n}\rightarrow
\infty $. Let $\hat{\theta}$ stand for any of the estimators $\hat{\theta}%
_{S}(\eta _{n})$, $\hat{\theta}_{H}(\eta _{n})$, or $\hat{\theta}_{A}(\eta
_{n})$. Let $d$ be a real number, and define the interval $D_{n}=[\hat{\theta%
}-d\eta _{n},\hat{\theta}+d\eta _{n}]$. If $d>1$, the interval $D_{n}$ has
infimal coverage probability converging to $1$, i.e., 
\begin{equation*}
\lim_{n\rightarrow \infty }\inf_{\theta \in \mathbb{R}}P_{n,\theta }(\theta
\in D_{n})=1\text{.}
\end{equation*}%
If $d<1$, 
\begin{equation*}
\lim_{n\rightarrow \infty }\inf_{\theta \in \mathbb{R}}P_{n,\theta }(\theta
\in D_{n})=0\text{.}
\end{equation*}
\end{proposition}

The asymptotic distributional results in the above proposition do not
provide information on the case $d=1$. However, from the finite-sample
results in Section \ref{inf_cov_prob} we see that in this case the infimal
coverage probability of $D_{n}$ converges to $1/2$.

Since the interval $D_{n}$ for $d>1$ has asymptotic infimal coverage
probability equal to one, one may wonder how much cruder this interval is
compared to the finite-sample intervals $C_{S,n}^{\ast }$, $C_{H,n}^{\ast }$%
, and $C_{A,n}^{\ast }$ constructed in Section \ref{symm}, which have
infimal coverage probability equal to a prespecified level $\delta $, $%
0<\delta <1$: The ratio of the half-length of $D_{n}$ to the half-length of
the corresponding interval $C_{S,n}^{\ast }$, $C_{H,n}^{\ast }$, and $%
C_{A,n}^{\ast }$ is 
\begin{equation*}
d(1+O(n^{-1/2}\eta _{n}^{-1}))=d(1+o(1))
\end{equation*}%
as can be seen from equations (\ref{asy_half-length_S}), (\ref%
{asy_half-length_H}), and (\ref{asy_half-length_A}). Since $d$ can be chosen
arbitrarily close to one, this ratio can be made arbitrarily close to one.
This may sound somewhat strange, since we are comparing an interval with
asymptotic infimal coverage probability $1$ with the shortest finite-sample
confidence intervals that have a fixed infimal coverage probability $\delta $
less than $1$. The reason for this phenomenon is that, in the relevant
moving-parameter asymptotic framework, the distribution of $\hat{\theta}%
-\theta $ is made up of a bias-component which in the worst case is of the
order $\eta _{n}$ and a random component which is of the order $n^{-1/2}$.
Since $\eta _{n}\rightarrow 0$ and $n^{1/2}\eta _{n}\rightarrow \infty $,
the deterministic bias-component dominates the random component. This can
also be gleaned from equations (\ref{asy_half-length_S}), (\ref%
{asy_half-length_H}), and (\ref{asy_half-length_A}), where the level $\delta 
$ enters the formula for the half-length only in the lower order term.

We note that using Theorem 19 in P\"{o}tscher and Leeb (2009) the same proof
immediately shows that Proposition \ref{asy} also holds for the smoothly
clipped absolute deviation (SCAD) estimator when tuned to possess the
`sparsity property'. In fact, the argument in the proof of the above
proposition can be applied to a large class of post-model-selection
estimators based on a consistent model selection procedure.

\begin{remark}
(i) Suppose $D_{n}^{\prime }=[\hat{\theta}-d_{1}\eta _{n},\hat{\theta}%
+d_{2}\eta _{n}]$ where $\hat{\theta}$ stands for any of the estimators $%
\hat{\theta}_{S}$, $\hat{\theta}_{H}$, or $\hat{\theta}_{A}$. If $\min
(d_{1},d_{2})>1$, then the limit of the infimal coverage probability of $%
D_{n}^{\prime }$ is $1$; if $\max (d_{1},d_{2})<1$ then this limit is zero.
This follows immediately from an inspection of the proof of Proposition \ref%
{asy}.

(ii) Proposition \ref{asy} also remains correct if $D_{n}$ is replaced by
the corresponding open interval. A similar comment applies to the open
version of $D_{n}^{\prime }$.
\end{remark}

\section{Confidence Intervals: Unknown-Variance Case}

In this section we consider the case where the variance $\sigma ^{2}$ is
unknown, $n\geq 2$, and we are interested in the coverage properties of
intervals of the form $[\tilde{\theta}-\hat{\sigma}a_{n},\tilde{\theta}+\hat{%
\sigma}a_{n}]$ where $a_{n}$ is a nonnegative real number and $\tilde{\theta}
$ stands for any one of the estimators $\tilde{\theta}_{H}=\tilde{\theta}%
_{H}(\eta _{n})$, $\tilde{\theta}_{S}=\tilde{\theta}_{S}(\eta _{n})$, or $%
\tilde{\theta}_{A}=\tilde{\theta}_{A}(\eta _{n})$. For brevity we only
consider symmetric intervals. A similar argument as in the known-variance
case shows that we can assume without loss of generality that $\sigma =1$,
and we shall do so in the sequel; in particular, this argument shows that
the infimum with respect to $\theta $ of the coverage probability does not
depend on $\sigma $.

\subsection{Soft-thresholding}

Consider the interval $E_{S,n}=\left[ \tilde{\theta}_{S}-\hat{\sigma}a_{n},%
\tilde{\theta}_{S}+\hat{\sigma}a_{n}\right] $ where $a_{n}$ is a nonnegative
real number and $\tilde{\theta}_{S}=\tilde{\theta}_{S}(\eta _{n})$. We then
have%
\begin{equation*}
P_{n,\theta }\left( \theta \in E_{S,n}\right) =\int_{0}^{\infty }P_{n,\theta
}\left( \theta \in E_{S,n}\left\vert \hat{\sigma}=s\right. \right) h_{n}(s)ds
\end{equation*}%
where $h_{n}$ is the density of $\hat{\sigma}$, i.e., $h_{n}$ is the density
of the square root of a chi-square distributed random variable with $n-1$
degrees of freedom divided by the degrees of freedom. In view of
independence of $\hat{\sigma}$ and $\bar{y}$ we obtain the following
representation of the finite-sample coverage probability%
\begin{eqnarray}
P_{n,\theta }\left( \theta \in E_{S,n}\right)  &=&\int_{0}^{\infty
}P_{n,\theta }\left( \theta \in \left[ \hat{\theta}_{S}(s\eta _{n})-sa_{n},%
\hat{\theta}_{S}(s\eta _{n})+sa_{n}\right] \right) h_{n}(s)ds  \notag \\
&=&\int_{0}^{\infty }p_{S,n}\left( \theta ;1,s\eta _{n},sa_{n},sa_{n}\right)
h_{n}(s)ds  \label{cov_unknown_S}
\end{eqnarray}%
where $p_{S,n}$ is given in (\ref{cov_S}) in the Appendix.

We next determine the infimal coverage probability of $E_{S,n}$ in finite
samples: It follows from (\ref{cov_S}), the dominated convergence theorem,
and symmetry of the standard normal distribution that%
\begin{eqnarray}
\inf_{\theta \in \mathbb{R}}P_{n,\theta }\left( \theta \in E_{S,n}\right)
&\leq &\lim_{\theta \rightarrow \infty }\int_{0}^{\infty }p_{S,n}\left(
\theta ;1,s\eta _{n},sa_{n},sa_{n}\right) h_{n}(s)ds  \notag \\
&=&\int_{0}^{\infty }\lim_{\theta \rightarrow \infty }p_{S,n}\left( \theta
;1,s\eta _{n},sa_{n},sa_{n}\right) h_{n}(s)ds  \notag \\
&=&\int_{0}^{\infty }[\Phi (n^{1/2}s(a_{n}-\eta _{n}))-\Phi
(n^{1/2}s(-a_{n}-\eta _{n}))]h_{n}(s)ds  \notag \\
&=&T_{n-1}(n^{1/2}(a_{n}-\eta _{n}))-T_{n-1}(n^{1/2}(-a_{n}-\eta _{n})),
\label{upper}
\end{eqnarray}%
where $T_{n-1}$ is the cdf of a Student $t$-distribution with $n-1$ degrees
of freedom. Furthermore, (\ref{infimal_soft_asym}) shows that 
\begin{equation*}
p_{S,n}\left( \theta ;1,s\eta _{n},sa_{n},sa_{n}\right) \geq \Phi
(n^{1/2}s(a_{n}-\eta _{n}))-\Phi (n^{1/2}s(-a_{n}-\eta _{n}))
\end{equation*}%
holds and whence we obtain from (\ref{cov_unknown_S}) and (\ref{upper}) the
following expression for the infimal coverage probability of $E_{S,n}$:%
\begin{equation}
\inf_{\theta \in \mathbb{R}}P_{n,\theta }\left( \theta \in E_{S,n}\right)
=T_{n-1}(n^{1/2}(a_{n}-\eta _{n}))-T_{n-1}(n^{1/2}(-a_{n}-\eta _{n}))
\label{finite}
\end{equation}%
for every $n\geq 2$. Remark \ref{open} shows that the same relation is true
for the corresponding open and half-open intervals.

Relation (\ref{finite}) shows the following: suppose $1/2\leq \delta <1$ and 
$a_{n,S}^{\ast }$ solves (\ref{short_a_S}), i.e., the corresponding interval 
$C_{S,n}^{\ast }$ has infimal coverage probability equal to $\delta $. Let $%
a_{n,S}^{\ast \ast }$ be the (unique)\ solution to 
\begin{equation*}
T_{n-1}(n^{1/2}(a_{n}-\eta _{n}))-T_{n-1}(n^{1/2}(-a_{n}-\eta _{n}))=\delta ,
\end{equation*}%
i.e., the corresponding interval $E_{S,n}^{\ast \ast }=\left[ \tilde{\theta}%
_{S}-\hat{\sigma}a_{n,S}^{\ast \ast },\tilde{\theta}_{S}+\hat{\sigma}%
a_{n,S}^{\ast \ast }\right] $ has infimal coverage probability equal to $%
\delta $. Then $a_{n,S}^{\ast \ast }\geq a_{n,S}^{\ast }$ holds in view of
Lemma \ref{l_5} in the Appendix. I.e., given the same infimal coverage
probability $\delta \geq 1/2$, the expected length of the interval $%
E_{S,n}^{\ast \ast }$ based on $\tilde{\theta}_{S}$ is not smaller than the
length of the interval $C_{S,n}^{\ast }$\ based on $\hat{\theta}_{S}$.

Since $\left\Vert \Phi -T_{n-1}\right\Vert _{\infty }=\sup_{x\in \mathbb{R}%
}\left\vert \Phi (x)-T_{n-1}(x)\right\vert \rightarrow 0$ for $n\rightarrow
\infty $ holds by Polya's theorem, the following result is an immediate
consequence of (\ref{finite}), Proposition \ref{lasso}, and Remark \ref{open}%
.

\begin{theorem}
\label{soft_unknown}For every sequence $a_{n}$ of nonnegative real numbers
we have with $E_{S,n}=\left[ \tilde{\theta}_{S}-\hat{\sigma}a_{n},\tilde{%
\theta}_{S}+\hat{\sigma}a_{n}\right] $ and $C_{S,n}=\left[ \hat{\theta}%
_{S}-a_{n},\hat{\theta}_{S}+a_{n}\right] $ that%
\begin{equation*}
\inf_{\theta \in \mathbb{R}}P_{n,\theta }\left( \theta \in E_{S,n}\right)
-\inf_{\theta \in \mathbb{R}}P_{n,\theta }\left( \theta \in C_{S,n}\right)
\rightarrow 0
\end{equation*}%
as $n\rightarrow \infty $. The analogous results hold for the corresponding
open and half-open intervals.
\end{theorem}

We discuss this theorem together with the parallel results for
hard-thresholding and adaptive LASSO based intervals in Section \ref{disc}.

\subsection{Hard-thresholding}

Consider the interval $E_{H,n}=\left[ \tilde{\theta}_{H}-\hat{\sigma}a_{n},%
\tilde{\theta}_{H}+\hat{\sigma}a_{n}\right] $ where $a_{n}$ is a nonnegative
real number and $\tilde{\theta}_{H}=\tilde{\theta}_{H}(\eta _{n})$. We then
have analogously as in the preceding subsection that%
\begin{equation*}
P_{n,\theta }\left( \theta \in E_{H,n}\right) =\int_{0}^{\infty
}p_{H,n}\left( \theta ;1,s\eta _{n},sa_{n},sa_{n}\right) h_{n}(s)ds.
\end{equation*}%
Note that $p_{H,n}\left( \theta ;1,s\eta _{n},sa_{n},sa_{n}\right) $ is
symmetric in $\theta $ and for $\theta \geq 0$ is given by (see P\"{o}tscher
(2009))%
\begin{eqnarray*}
&&p_{H,n}\left( \theta ;1,s\eta _{n},sa_{n},sa_{n}\right) \\
&=&\left\{ \Phi (n^{1/2}(-\theta +s\eta _{n}))-\Phi (n^{1/2}(-\theta -s\eta
_{n}))\right\} \boldsymbol{1}\left( 0\leq \theta \leq sa_{n}\right) \\
&&+\max \left[ 0,\Phi (n^{1/2}sa_{n})-\Phi (n^{1/2}(-\theta +s\eta _{n}))%
\right] \boldsymbol{1}\left( sa_{n}<\theta \leq s\eta _{n}+sa_{n}\right) \\
&&+\left\{ \Phi (n^{1/2}sa_{n})-\Phi (-n^{1/2}sa_{n})\right\} \boldsymbol{1}%
\left( s\eta _{n}+sa_{n}<\theta \right)
\end{eqnarray*}%
in case $\eta _{n}>2a_{n}$, by%
\begin{eqnarray*}
&&p_{H,n}\left( \theta ;1,s\eta _{n},sa_{n},sa_{n}\right) \\
&=&\left\{ \Phi (n^{1/2}(-\theta +s\eta _{n})-\Phi (n^{1/2}(-\theta -s\eta
_{n}))\right\} \boldsymbol{1}\left( 0\leq \theta \leq s\eta
_{n}-sa_{n}\right) \\
&&+\left\{ \Phi (n^{1/2}sa_{n})-\Phi (n^{1/2}(-\theta -s\eta _{n}))\right\} 
\boldsymbol{1}\left( s\eta _{n}-sa_{n}<\theta \leq sa_{n}\right) \\
&&+\left\{ \Phi (n^{1/2}sa_{n})-\Phi (n^{1/2}(-\theta +s\eta _{n}))\right\} 
\boldsymbol{1}\left( sa_{n}<\theta \leq s\eta _{n}+sa_{n}\right) \\
&&+\left\{ \Phi (n^{1/2}sa_{n})-\Phi (-n^{1/2}sa_{n})\right\} \boldsymbol{1}%
\left( s\eta _{n}+sa_{n}<\theta \right)
\end{eqnarray*}%
if $a_{n}\leq \eta _{n}\leq 2a_{n}$, and by 
\begin{eqnarray*}
&&p_{H,n}\left( \theta ;1,s\eta _{n},sa_{n},sa_{n}\right) \\
&=&\left\{ \Phi (n^{1/2}sa_{n})-\Phi (-n^{1/2}sa_{n})\right\} \left\{ 
\boldsymbol{1}\left( 0\leq \theta \leq sa_{n}-s\eta _{n}\right) +\boldsymbol{%
1}\left( s\eta _{n}+sa_{n}<\theta \right) \right\} \\
&&+\left\{ \Phi (n^{1/2}sa_{n})-\Phi (n^{1/2}(-\theta -s\eta _{n}))\right\} 
\boldsymbol{1}\left( sa_{n}-s\eta _{n}<\theta \leq sa_{n}\right) \\
&&+\left\{ \Phi (n^{1/2}sa_{n})-\Phi (n^{1/2}(-\theta +s\eta _{n}))\right\} 
\boldsymbol{1}\left( sa_{n}<\theta \leq s\eta _{n}+sa_{n}\right)
\end{eqnarray*}%
if $\eta _{n}<a_{n}$. In the subsequent theorems we consider only the case
where $\eta _{n}\rightarrow 0$ as this is the only interesting case from an
asymptotic perspective: note that any of the penalized maximum likelihood
estimators considered in this paper is inconsistent for $\theta $ if $\eta
_{n}$ does not converge to zero.

\begin{theorem}
\label{hard_unknown}Suppose $\eta _{n}\rightarrow 0$. For every sequence $%
a_{n}$ of nonnegative real numbers we have with $E_{H,n}=\left[ \tilde{\theta%
}_{H}-\hat{\sigma}a_{n},\tilde{\theta}_{H}+\hat{\sigma}a_{n}\right] $ and $%
C_{H,n}=\left[ \hat{\theta}_{H}-a_{n},\hat{\theta}_{H}+a_{n}\right] $ that%
\begin{equation*}
\inf_{\theta \in \mathbb{R}}P_{n,\theta }\left( \theta \in E_{H,n}\right)
-\inf_{\theta \in \mathbb{R}}P_{n,\theta }\left( \theta \in C_{H,n}\right)
\rightarrow 0
\end{equation*}%
as $n\rightarrow \infty $. The analogous results hold for the corresponding
open and half-open intervals.
\end{theorem}

\subsection{Adaptive LASSO}

Consider the interval $E_{A,n}=[\tilde{\theta}_{A}-\hat{\sigma}a_{n},\tilde{%
\theta}_{A}+\hat{\sigma}a_{n}]$ where $a_{n}$ is a nonnegative real number
and $\tilde{\theta}_{A}=\tilde{\theta}_{A}(\eta _{n})$. We then have
analogously as in the preceding subsections that 
\begin{equation*}
P_{n,\theta }(\theta \in E_{A,n})=\int_{0}^{\infty }p_{A,n}(\theta ;1,s\eta
_{n},sa_{n},sa_{n})h_{n}(s)ds
\end{equation*}%
where $p_{A,n}$ is given in (\ref{cov}) in the Appendix.

\begin{theorem}
\label{adaptive_unknown}Suppose $\eta _{n}\rightarrow 0$. For every sequence 
$a_{n}$ of nonnegative real numbers we have with $E_{A,n}=\left[ \tilde{%
\theta}_{A}-\hat{\sigma}a_{n},\tilde{\theta}_{A}+\hat{\sigma}a_{n}\right] $
and $C_{A,n}=\left[ \hat{\theta}_{A}-a_{n},\hat{\theta}_{A}+a_{n}\right] $
that%
\begin{equation*}
\inf_{\theta \in \mathbb{R}}P_{n,\theta }\left( \theta \in E_{A,n}\right)
-\inf_{\theta \in \mathbb{R}}P_{n,\theta }\left( \theta \in C_{A,n}\right)
\rightarrow 0
\end{equation*}%
as $n\rightarrow \infty $. The analogous results hold for the corresponding
open and half-open intervals.
\end{theorem}

\subsection{Discussion\label{disc}}

Theorems \ref{soft_unknown}, \ref{hard_unknown}, and \ref{adaptive_unknown}
show that the results in Section \ref{finite_sample} carry over to the
unknown-variance case in an asymptotic sense: For example, suppose $0<\delta
<1$, and $a_{n,S}$ ($a_{n,H}$, $a_{n,A}$, respectively) is such that $%
E_{S,n} $ ($E_{H,n}$, $E_{A,n}$, respectively) has infimal coverage
probability converging to $\delta $. Then, for a regime where $n^{1/2}\eta
_{n}\rightarrow e$ with $0\leq e<\infty $, it follows that $n^{1/2}a_{n,S}$, 
$n^{1/2}a_{n,H}$, and $n^{1/2}a_{n,A}$ have limits that solve (\ref%
{asy_half-lenght_S_0})-(\ref{asy_half-lenght_A_0}), respectively; that is,
they have the same limits as $n^{1/2}a_{n,S}^{\ast }$, $n^{1/2}a_{n,H}^{\ast
}$, and $n^{1/2}a_{n,A}^{\ast }$, which are $n^{1/2}$ times the half-length
of the shortest $\delta $-confidence intervals $C_{S,n}^{\ast }$, $%
C_{H,n}^{\ast }$, and $C_{A,n}^{\ast }$, respectively, in the known-variance
case. Furthermore, for a regime where $n^{1/2}\eta _{n}\rightarrow \infty $
it follows that $a_{n,S}$, $a_{n,H}$, and $a_{n,A}$ satisfy (\ref%
{asy_half-length_S})-(\ref{asy_half-length_A}), respectively (where we also
assume $\eta _{n}\rightarrow 0$ for hard-thresholding and the adaptive
LASSO). Hence, $a_{n,S}$, $a_{n,H}$, and $a_{n,A}$ on the one hand, and $%
a_{n,S}^{\ast }$, $a_{n,H}^{\ast }$, and $a_{n,A}^{\ast }$ on the other hand
have again the same asymptotic behavior. Furthermore, Theorems \ref%
{soft_unknown}, \ref{hard_unknown}, and \ref{adaptive_unknown} show that
Proposition \ref{asy} immediately carries over to the unknown-variance case.

\appendix{}

\section{Appendix}

\textbf{Proof of Proposition \ref{lasso}: }Using the expression for the
finite sample distribution of $n^{1/2}(\hat{\theta}_{S}-\theta )$ given in P%
\"{o}tscher and Leeb (2009) and noting that this distribution function has a
jump at $-n^{1/2}\theta $ we obtain%
\begin{eqnarray}
p_{S,n}(\theta ) &=&[\Phi (n^{1/2}(a_{n}-\eta _{n}))-\Phi
(n^{1/2}(-b_{n}-\eta _{n}))]\boldsymbol{1}(\theta <-a_{n})  \notag \\
&+&[\Phi (n^{1/2}(a_{n}+\eta _{n}))-\Phi (n^{1/2}(-b_{n}-\eta _{n}))]%
\boldsymbol{1}(-a_{n}\leq \theta \leq b_{n})  \notag \\
&+&[\Phi (n^{1/2}(a_{n}+\eta _{n}))-\Phi (n^{1/2}(-b_{n}+\eta _{n}))]%
\boldsymbol{1}(b_{n}<\theta ).  \label{cov_S}
\end{eqnarray}%
It follows that $\inf_{\theta \in \mathbb{R}}p_{S,n}(\theta )$ is as given
in the proposition. $\ \blacksquare $

\bigskip

\noindent \textbf{Proof of Proposition \ref{adLASSOinf}: }The distribution
function $F_{A,n,\theta }(x)=P_{n,\theta }(n^{1/2}(\hat{\theta}_{A}-\theta
)\leq x)$ of the adaptive LASSO estimator is given by 
\begin{eqnarray*}
\boldsymbol{1}(x+n^{1/2}\theta  &\geq &0)\Phi \left( -((n^{1/2}\theta -x)/2)+%
\sqrt{((n^{1/2}\theta +x)/2)^{2}+n\eta _{n}^{2}}\right) + \\
\boldsymbol{1}(x+n^{1/2}\theta  &<&0)\Phi \left( -((n^{1/2}\theta -x)/2)-%
\sqrt{((n^{1/2}\theta +x)/2)^{2}+n\eta _{n}^{2}}\right) 
\end{eqnarray*}%
(see P\"{o}tscher and Schneider (2009)). Hence, the coverage probability $%
p_{A,n}(\theta )=F_{A,n,\theta }(n^{1/2}a_{n})-\lim_{x\rightarrow
(-n^{1/2}b_{n})_{-}}F_{A,n,\theta }(x)$ is%
\begin{equation}
p_{A,n}(\theta )=\left\{ 
\begin{array}{ll}
\Phi \left( n^{1/2}\gamma ^{(-)}(\theta ,-a_{n})\right) -\Phi \left(
n^{1/2}\gamma ^{(-)}(\theta ,b_{n})\right)  & \text{if }\;\theta <-a_{n} \\ 
\Phi \left( n^{1/2}\gamma ^{(+)}(\theta ,-a_{n})\right) -\Phi \left(
n^{1/2}\gamma ^{(-)}(\theta ,b_{n})\right)  & \text{if }\;-a_{n}\leq \theta
\leq b_{n} \\ 
\Phi \left( n^{1/2}\gamma ^{(+)}(\theta ,-a_{n})\right) -\Phi \left(
n^{1/2}\gamma ^{(+)}(\theta ,b_{n})\right)  & \text{if }\;\theta >b_{n}.%
\end{array}%
\right.   \label{cov}
\end{equation}%
Here%
\begin{eqnarray}
\gamma ^{(-)}(\theta ,x) &=&-((\theta +x)/2)-\sqrt{((\theta -x)/2)^{2}+\eta
_{n}^{2}}  \label{gamma1} \\
\gamma ^{(+)}(\theta ,x) &=&-((\theta +x)/2)+\sqrt{((\theta -x)/2)^{2}+\eta
_{n}^{2}},  \label{gamma2}
\end{eqnarray}%
which are clearly smooth functions of $(\theta ,x)$. Observe that $\gamma
^{(-)}$ and $\gamma ^{(+)}$ are nonincreasing in $\theta \in \mathbb{R}$
(for every $x\in \mathbb{R}$). As a consequence, we obtain for $-a_{n}\leq
\theta \leq b_{n}$ the lower bound%
\begin{eqnarray}
p_{A,n}(\theta ) &\geq &\Phi \left( n^{1/2}\gamma
^{(+)}(b_{n},-a_{n})\right) -\Phi \left( n^{1/2}\gamma
^{(-)}(-a_{n},b_{n})\right)   \notag \\
&=&\Phi \left( n^{1/2}\left( (a_{n}-b_{n})/2+\sqrt{((a_{n}+b_{n})/2)^{2}+%
\eta _{n}^{2}}\right) \right)   \notag \\
&&-\Phi \left( n^{1/2}\left( (a_{n}-b_{n})/2-\sqrt{((a_{n}+b_{n})/2)^{2}+%
\eta _{n}^{2}}\right) \right) .  \label{lower_bound}
\end{eqnarray}

Consider first the case where $a_{n}\leq b_{n}$. We then show that $%
p_{A,n}(\theta )$ is nonincreasing on $(-\infty ,-a_{n})$: The derivative $%
dp_{A,n}(\theta )/d\theta $ is given by 
\begin{eqnarray*}
&&dp_{A,n}(\theta )/d\theta = \\
&&n^{1/2}[\phi (n^{1/2}\gamma ^{(-)}(\theta ,-a_{n}))\partial \gamma
^{(-)}(\theta ,-a_{n})/\partial \theta -\phi (n^{1/2}\gamma ^{(-)}(\theta
,b_{n}))\partial \gamma ^{(-)}(\theta ,b_{n})/\partial \theta ]
\end{eqnarray*}%
where $\phi $ denotes the standard normal density function. Using the
relation $a_{n}\leq b_{n}$, elementary calculations show that 
\begin{equation*}
\partial \gamma ^{(-)}(\theta ,-a_{n})/\partial \theta \leq \partial \gamma
^{(-)}(\theta ,b_{n})/\partial \theta \text{ \ \ \ for }\theta \in (-\infty
,-a_{n})\text{.}
\end{equation*}%
Furthermore, given $a_{n}\leq b_{n}$, it is not too difficult to see that $%
\left\vert \gamma ^{(-)}(\theta ,-a_{n})\right\vert \leq \left\vert \gamma
^{(-)}(\theta ,b_{n})\right\vert $ for $\theta \in (-\infty ,-a_{n})$ (cf.
Lemma \ref{l_1} below), which implies that%
\begin{equation*}
\phi (n^{1/2}\gamma ^{(-)}(\theta ,-a_{n}))\geq \phi (n^{1/2}\gamma
^{(-)}(\theta ,b_{n})).
\end{equation*}%
The last two displays together with the fact that $\partial \gamma
^{(-)}(\theta ,-a_{n})/\partial \theta $ as well as $\partial \gamma
^{(-)}(\theta ,b_{n})/\partial \theta $ are less than or equal to zero,
imply that $dp_{A,n}(\theta )/d\theta \leq 0$ on $(-\infty ,-a_{n})$. This
proves that 
\begin{equation*}
\inf_{\theta <-a_{n}}p_{A,n}(\theta )=\lim_{\theta \rightarrow
(-a_{n})_{-}}p_{A,n}(\theta )=c
\end{equation*}%
with 
\begin{equation}
c=\Phi \left( n^{1/2}(a_{n}-\eta _{n})\right) -\Phi \left( n^{1/2}\left(
(a_{n}-b_{n})/2-\sqrt{((a_{n}+b_{n})/2)^{2}+\eta _{n}^{2}}\right) \right) .
\label{c}
\end{equation}%
Since the lower bound given in (\ref{lower_bound}) is not less than $c$, we
have 
\begin{equation*}
\inf_{\theta \leq b_{n}}p_{A,n}(\theta )=\inf_{\theta <-a_{n}}p_{A,n}(\theta
)=c.
\end{equation*}%
It remains to show that $p_{A,n}(\theta )\geq c$ for $\theta >b_{n}$. From (%
\ref{cov}) and (\ref{c}) after rearranging terms we obtain for $\theta
>b_{n} $%
\begin{eqnarray*}
p_{A,n}(\theta )-c &=&\left[ \Phi (n^{1/2}\gamma ^{(+)}(\theta
,-a_{n}))-\Phi (n^{1/2}\gamma ^{(-)}(-a_{n},-a_{n}))\right] - \\
&&\left[ \Phi (n^{1/2}\gamma ^{(+)}(\theta ,b_{n}))-\Phi (n^{1/2}\gamma
^{(-)}(-a_{n},b_{n}))\right] .
\end{eqnarray*}%
It is elementary to show that $\gamma ^{(+)}(\theta ,-a_{n}))\geq \gamma
^{(-)}(-a_{n},-a_{n})=a_{n}-\eta _{n}$ and $\gamma ^{(+)}(\theta
,b_{n}))\geq \gamma ^{(-)}(-a_{n},b_{n})$. We next show that%
\begin{equation}
\gamma ^{(+)}(\theta ,-a_{n})-\gamma ^{(-)}(-a_{n},-a_{n}))\geq \gamma
^{(+)}(\theta ,b_{n})-\gamma ^{(-)}(-a_{n},b_{n}).  \label{comp_lenght}
\end{equation}%
To establish this note that (\ref{comp_lenght}) can equivalently be
rewritten as%
\begin{equation}
f(0)+f((\theta +a_{n})/2)\geq f((\theta -b_{n})/2)+f((a_{n}+b_{n})/2)
\label{ineq_f}
\end{equation}%
where $f(x)=(x^{2}+\eta _{n}^{2})^{1/2}$. Observe that $0\leq (\theta
-b_{n})/2\leq (\theta +a_{n})/2$ holds since $0\leq a_{n}\leq b_{n}<\theta $%
. Writing $(\theta -b_{n})/2$ as $\lambda (\theta +a_{n})/2+(1-\lambda )0$
with $0\leq \lambda \leq 1$ gives $(a_{n}+b_{n})/2=(1-\lambda )(\theta
+a_{n})/2+\lambda 0$. Because $f$ is convex, the inequality (\ref{ineq_f})
and hence (\ref{comp_lenght}) follows.

Next observe that in case $a_{n}\geq \eta _{n}$ we have (using monotonicity
of $\gamma ^{(+)}(\theta ,b_{n})$)%
\begin{equation}
0\leq \gamma ^{(-)}(-a_{n},-a_{n}))=a_{n}-\eta _{n}\leq b_{n}-\eta
_{n}=-\gamma ^{(+)}(b_{n},b_{n})\leq -\gamma ^{(+)}(\theta ,b_{n})
\label{ineq_1}
\end{equation}%
for $\theta >b_{n}$. In case $a_{n}<\eta _{n}$ we have (using $\gamma
^{(-)}(\theta ,x)=\gamma ^{(-)}(x,\theta )$ and monotonicity of $\gamma
^{(-)}$ in its first argument)%
\begin{equation}
\gamma ^{(-)}(-a_{n},b_{n})\leq \gamma ^{(-)}(-a_{n},-a_{n})=a_{n}-\eta
_{n}<0,  \label{ineq_3}
\end{equation}%
and (using monotonicity of $\gamma ^{(+)}$) 
\begin{equation}
\gamma ^{(-)}(-a_{n},b_{n})\leq -\gamma ^{(+)}(b_{n},-a_{n})\leq -\gamma
^{(+)}(\theta ,-a_{n})  \label{ineq_4}
\end{equation}%
for $\theta >b_{n}$. Applying Lemma \ref{l_2} below with $\alpha
=n^{1/2}\gamma ^{(-)}(-a_{n},-a_{n})$, $\beta =n^{1/2}\gamma ^{(+)}(\theta
,-a_{n})$, $\gamma =n^{1/2}\gamma ^{(-)}(-a_{n},b_{n})$, and $\delta
=n^{1/2}\gamma ^{(+)}(\theta ,b_{n})$ and using (\ref{comp_lenght})-(\ref%
{ineq_4}), establishes $p_{A,n}(\theta )-c\geq 0$. This completes the proof
in case $a_{n}\leq b_{n}$.

The case $a_{n}>b_{n}$ follows from the observation that (\ref{cov}) remains
unchanged if $a_{n}$ and $b_{n}$ are interchanged and $\theta $ is replaced
by $-\theta $. $\ \blacksquare $

\begin{lemma}
\label{l_1} Suppose $a_{n}\leq b_{n}$. Then $\left\vert \gamma ^{(-)}(\theta
,-a_{n})\right\vert \leq \left\vert \gamma ^{(-)}(\theta ,b_{n})\right\vert $
holds for $\theta \in (-\infty ,-a_{n})$.
\end{lemma}

\begin{proof}
Squaring both sides of the claimed inequality shows that the claim is
equivalent to%
\begin{equation*}
a_{n}^{2}/2-(a_{n}-\theta )\sqrt{((a_{n}+\theta )/2)^{2}+\eta ^{2}}\leq
b_{n}^{2}/2+(b_{n}+\theta )\sqrt{((b_{n}-\theta )/2)^{2}+\eta ^{2}}.
\end{equation*}%
But, for $\theta <-a_{n}$, the left-hand side of the preceding display is
not larger than%
\begin{equation*}
a_{n}^{2}/2+(a_{n}+\theta )\sqrt{((a_{n}-\theta )/2)^{2}+\eta ^{2}}.
\end{equation*}%
Since $a_{n}^{2}/2\leq b_{n}^{2}/2$, it hence suffices to show that 
\begin{equation*}
-(a_{n}+\theta )\sqrt{((a_{n}-\theta )/2)^{2}+\eta ^{2}}\geq -(b_{n}+\theta )%
\sqrt{((b_{n}-\theta )/2)^{2}+\eta ^{2}}
\end{equation*}%
for $\theta <-a_{n}$. This is immediately seen by distinguishing the cases
where $-b_{n}\leq \theta <-a_{n}$ and where $\theta <-b_{n}$, and observing
that $a_{n}\leq b_{n}$.
\end{proof}

\bigskip

The following lemma is elementary to prove.

\begin{lemma}
\label{l_2} Suppose $\alpha $, $\beta $, $\gamma $, and $\delta $ are real
numbers satisfying $\alpha \leq \beta $, $\gamma \leq \delta $, and $\beta
-\alpha \geq \delta -\gamma $. If $0\leq \alpha \leq -\delta $, or if $%
\gamma \leq \alpha \leq 0$ and $\gamma \leq -\beta $, then $\Phi (\beta
)-\Phi (\alpha )\geq \Phi (\delta )-\Phi (\gamma )$.
\end{lemma}

\bigskip

\noindent \textbf{Proof of Theorem \ref{short}: }(a) Since $\delta $ is
positive, any solution to (\ref{short_a_S}) has to be positive. Now the
equation (\ref{short_a_S}) has a unique solution $a_{n,S}^{\ast }$, since (%
\ref{short_a_S}) as a function of $a_{n}\in \lbrack 0,\infty )$ is easily
seen to be strictly increasing with range $[0,1)$. Furthermore, the infimal
coverage probability (\ref{infimal_soft_asym}) is a continuous function of
the pair $(a_{n},b_{n})$ on $[0,\infty )\times \lbrack 0,\infty )$. Let $%
K\subseteq \lbrack 0,\infty )\times \lbrack 0,\infty )$ consist of all pairs 
$(a_{n},b_{n})$ such that (i) the corresponding interval $[\hat{\theta}%
_{S}-a_{n},\hat{\theta}_{S}+b_{n}]$ has infimal coverage probability not
less than $\delta $, and (ii) the length $a_{n}+b_{n}$ is less than or equal 
$2a_{n,S}^{\ast }$. Then $K$ is compact. It is also nonempty as the pair $%
(a_{n,S}^{\ast },a_{n,S}^{\ast })$ belongs to $K$. Since the length $%
a_{n}+b_{n}$ is obviously continuous, it follows that there is a pair $%
(a_{n}^{o},b_{n}^{o})\in K$ having minimal length within $K$. Since
confidence sets corresponding to pairs not belonging to $K$ always have
length larger than $2a_{n,S}^{\ast }$, the pair $(a_{n}^{o},b_{n}^{o})$
gives rise to an interval with shortest length within the set of all
intervals with infimal coverage probability not less than $\delta $. We next
show that $a_{n}^{o}=b_{n}^{o}$ must hold: Suppose not, then we may assume
without loss of generality that $a_{n}^{o}<b_{n}^{o}$, since (\ref%
{infimal_soft_asym}) remains invariant under permutation of $a_{n}^{o}$ and $%
b_{n}^{o}$. But now increasing $a_{n}^{o}$ by $\varepsilon >0$ and
decreasing $b_{n}^{o}$ by the same amount such that $a_{n}^{o}+\varepsilon
<b_{n}^{o}-\varepsilon $ holds, will result in an interval of the same
length with infimal coverage probability%
\begin{equation*}
\Phi (n^{1/2}(a_{n}^{o}+\varepsilon -\eta _{n}))-\Phi
(n^{1/2}(-(b_{n}^{o}-\varepsilon )-\eta _{n})).
\end{equation*}%
This infimal coverage probability will be strictly larger than%
\begin{equation*}
\Phi (n^{1/2}(a_{n}^{o}-\eta _{n}))-\Phi (n^{1/2}(-b_{n}^{o}-\eta _{n}))\geq
\delta
\end{equation*}%
provided $\varepsilon $ is chosen sufficiently small. But then, by
continuity of the infimal coverage probability as a function of $a_{n}$ and $%
b_{n}$, the interval $[\hat{\theta}_{S}-a_{n}^{o}-\varepsilon ,\hat{\theta}%
_{S}+b_{n}^{\prime }-\varepsilon ]$ with $\varepsilon <b_{n}^{\prime
}<b_{n}^{o}$ will still have infimal coverage probability not less than $%
\delta $ as long as $b_{n}^{\prime }$ is sufficiently close to $b_{n}^{o}$;
at the same time this interval will be shorter than the interval $[\hat{%
\theta}_{S}-a_{n}^{o},\hat{\theta}_{S}+b_{n}^{o}]$. This leads to a
contradiction and establishes $a_{n}^{o}=b_{n}^{o}$. By what was said at the
beginning of the proof, it is now obvious that $%
a_{n}^{o}=b_{n}^{o}=a_{n,S}^{\ast }$ must hold, thus also establishing
uniqueness. The last claim is obvious in view of the construction of $%
a_{n,S}^{\ast }$.

(b) Since $\delta $ is positive, any solution to (\ref{short_a_H}) has to be
larger than $\eta _{n}/2$. Now equation (\ref{short_a_H}) has a unique
solution $a_{n,H}^{\ast }$, since (\ref{short_a_H}) as a function of $%
a_{n}\in \lbrack \eta _{n}/2,\infty )$ is easily seen to be strictly
increasing with range $[0,1)$. Furthermore, define $K$ similarly as in the
proof of part (a). Then, by the same reasoning as in (a), the set $K$ is
compact and non-empty, leading to a pair $(a_{n}^{o},b_{n}^{o})$ that gives
rise to an interval with shortest length within the set of all intervals
with infimal coverage probability not less than $\delta $. We next show that 
$a_{n}^{o}=b_{n}^{o}$ must hold: Suppose not, then we may again assume
without loss of generality that $a_{n}^{o}<b_{n}^{o}$. Note that $%
a_{n}^{o}+b_{n}^{o}>\eta _{n}$ must hold, since the infimal coverage
probability of the corresponding interval is positive by construction. Since
all this entails $\left\vert a_{n}^{o}-\eta _{n}\right\vert <b_{n}^{o}$,
increasing $a_{n}^{o}$ by $\varepsilon >0$ and decreasing $b_{n}^{o}$ by the
same amount such that $a_{n}^{o}+\varepsilon <b_{n}^{o}-\varepsilon $ holds,
will result in an interval of the same length with infimal coverage
probability%
\begin{eqnarray*}
\Phi (n^{1/2}(a_{n}^{o}+\varepsilon -\eta _{n}))-\Phi
(-n^{1/2}(b_{n}^{o}-\varepsilon )) &>& \\
\Phi (n^{1/2}(a_{n}^{o}-\eta _{n}))-\Phi (-n^{1/2}b_{n}^{o}) &\geq &\delta
\end{eqnarray*}%
provided $\varepsilon $ is chosen sufficiently small. By continuity of the
infimal coverage probability as a function of $a_{n}$ and $b_{n}$, the
interval $[\hat{\theta}_{S}-a_{n}^{o}-\varepsilon ,\hat{\theta}%
_{S}+b_{n}^{\prime }-\varepsilon ]$ with $\varepsilon <b_{n}^{\prime
}<b_{n}^{o}$ will still have infimal coverage probability not less than $%
\delta $ as long as $b_{n}^{\prime }$ is sufficiently close to $b_{n}^{o}$;
at the same time this interval will be shorter than the interval $[\hat{%
\theta}_{S}-a_{n}^{o},\hat{\theta}_{S}+b_{n}^{o}]$, leading to a
contradiction thus establishing $a_{n}^{o}=b_{n}^{o}$. As in (a) it now
follows that $a_{n}^{o}=b_{n}^{o}=a_{n,H}^{\ast }$ must hold, thus also
establishing uniqueness. The last claim is then obvious in view of the
construction of $a_{n,H}^{\ast }$.

(c) Since $\delta $ is positive, it is easy to see that any solution to (\ref%
{short_a_A}) has to be positive. Now equation (\ref{short_a_A}) has a unique
solution $a_{n,A}^{\ast }$, since (\ref{short_a_A}) as a function of $%
a_{n}\in \lbrack 0,\infty )$ is strictly increasing with range $[0,1)$.
Furthermore, the infimal coverage probability as given in Proposition \ref%
{adLASSOinf} is a continuous function of the pair $(a_{n},b_{n})$ on $%
[0,\infty )\times \lbrack 0,\infty )$. Define $K$ similarly as in the proof
of part (a). Then by the same reasoning as in (a), the set $K$ is compact
and non-empty, leading to a pair $(a_{n}^{o},b_{n}^{o})$ that gives rise to
an interval with shortest length within the set of all intervals with
infimal coverage probability not less than $\delta $. We next show that $%
a_{n}^{o}=b_{n}^{o}$ must hold: Suppose not, then we may again assume
without loss of generality that $a_{n}^{o}<b_{n}^{o}$. But now increasing $%
a_{n}^{o}$ by $\varepsilon >0$ and decreasing $b_{n}^{o}$ by the same amount
such that $a_{n}^{o}+\varepsilon <b_{n}^{o}-\varepsilon $ holds, will result
in an interval of the same length with infimal coverage probability%
\begin{eqnarray*}
&&\Phi (n^{1/2}(a_{n}^{o}+\varepsilon -\eta _{n}))-\Phi \left( n^{1/2}\left(
\varepsilon +(a_{n}^{o}-b_{n}^{o})/2-\sqrt{((a_{n}^{o}+b_{n}^{o})/2)^{2}+%
\eta _{n}^{2}}\right) \right) > \\
&&\Phi (n^{1/2}(a_{n}^{o}-\eta _{n}))-\Phi \left( n^{1/2}\left(
(a_{n}^{o}-b_{n}^{o})/2-\sqrt{((a_{n}^{o}+b_{n}^{o})/2)^{2}+\eta _{n}^{2}}%
\right) \right) \geq \delta ,
\end{eqnarray*}%
provided $\varepsilon $ is chosen sufficiently small. This is so since $%
a_{n}^{o}<b_{n}^{o}$ implies 
\begin{equation*}
\left\vert a_{n}^{o}-\eta _{n}\right\vert <\left\vert
(a_{n}^{o}-b_{n}^{o})/2-\sqrt{((a_{n}^{o}+b_{n}^{o})/2)^{2}+\eta _{n}^{2}}%
\right\vert 
\end{equation*}%
as is easily seen. But then, by continuity of the infimal coverage
probability as a function of $a_{n}$ and $b_{n}$, the interval $[\hat{\theta}%
_{S}-a_{n}^{o}-\varepsilon ,\hat{\theta}_{S}+b_{n}^{\prime }-\varepsilon ]$
with $\varepsilon <b_{n}^{\prime }<b_{n}^{o}$ will still have infimal
coverage probability not less than $\delta $ as long as $b_{n}^{\prime }$ is
sufficiently close to $b_{n}^{o}$; at the same time this interval will be
shorter than the interval $[\hat{\theta}_{S}-a_{n}^{o},\hat{\theta}%
_{S}+b_{n}^{o}]$. This leads to a contradiction and establishes $%
a_{n}^{o}=b_{n}^{o}$. As in (a) it now follows that $%
a_{n}^{o}=b_{n}^{o}=a_{n,A}^{\ast }$ must hold, thus also establishing
uniqueness. The last claim is obvious in view of the construction of $%
a_{n,A}^{\ast }$. $\ \blacksquare $

\bigskip

\noindent \textbf{Proof of Proposition \ref{asy}: }Let 
\begin{equation*}
c=\liminf_{n\rightarrow \infty }\inf_{\theta \in \mathbb{R}}P_{n,\theta
}\left( -d\leq \eta _{n}^{-1}(\hat{\theta}-\theta )\leq d\right) .
\end{equation*}%
By definition of $c$, we can find a subsequence $n_{k}$ and elements $\theta
_{n_{k}}\in \mathbb{R}$ such that%
\begin{equation*}
P_{n_{k},\theta _{n_{k}}}\left( -d\leq \eta _{n_{k}}^{-1}(\hat{\theta}%
-\theta _{n_{k}})\leq d\right) \rightarrow c
\end{equation*}%
for $k\rightarrow \infty $. Now, by Theorem 17 (for $\hat{\theta}=\hat{\theta%
}_{H}$), Theorem 18 (for $\hat{\theta}=\hat{\theta}_{S}$), and Remark 12 in P%
\"{o}tscher and Leeb (2009), and by Theorem 6 (for $\hat{\theta}=\hat{\theta}%
_{A}$) and Remark 7 in P\"{o}tscher and Schneider (2009), any accumulation
point of the distribution of $\eta _{n_{k}}^{-1}(\hat{\theta}-\theta
_{n_{k}})$ with respect to weak convergence is a probability measure
concentrated on $[-1,1]$. Since $d>1$, it follows that $c=1$ must hold,
which proves the first claim. We next prove the second claim. In view of
Theorem 17 (for $\hat{\theta}=\hat{\theta}_{H}$) and Theorem 18 (for $\hat{%
\theta}=\hat{\theta}_{S}$) in P\"{o}tscher and Leeb (2009), and in view of
Theorem 6 (for $\hat{\theta}=\hat{\theta}_{A}$) in P\"{o}tscher and
Schneider (2009) it is possible to choose a sequence $\theta _{n}\in \mathbb{%
R}$ such that the distribution of $\eta _{n}^{-1}(\hat{\theta}-\theta _{n})$
converges to point mass located at one of the endpoints of the interval $%
[-1,1]$. But then clearly 
\begin{equation*}
P_{n,\theta _{n}}\left( -d\leq \eta _{n}^{-1}(\hat{\theta}-\theta _{n})\leq
d\right) \rightarrow 0
\end{equation*}%
for $d<1$ which implies the second claim. $\ \blacksquare $

\bigskip

\noindent \textbf{Proof of Theorem \ref{hard_unknown}: }We prove the result
for the closed interval. Inspection of the proof together with Remark \ref%
{open}\ then gives the result for the open and half-open intervals.

Step 1: Observe that for every $s>0$ and $n\geq 2$ we have from the above
formulae for $p_{H,n}$ that%
\begin{equation*}
\lim_{\theta \rightarrow \infty }p_{H,n}\left( \theta ;1,s\eta
_{n},sa_{n},sa_{n}\right) =\Phi (n^{1/2}sa_{n})-\Phi (-n^{1/2}sa_{n}).
\end{equation*}%
By the dominated convergence theorem it follows that for $\theta \rightarrow
\infty $%
\begin{eqnarray*}
P_{n,\theta }\left( \theta \in E_{H,n}\right) &=&\int_{0}^{\infty
}p_{H,n}\left( \theta ;1,s\eta _{n},sa_{n},sa_{n}\right) h_{n}(s)ds \\
&\rightarrow &\int_{0}^{\infty }\left[ \Phi (n^{1/2}sa_{n})-\Phi
(-n^{1/2}sa_{n})\right] h_{n}(s)ds \\
&=&T_{n-1}(n^{1/2}a_{n})-T_{n-1}(-n^{1/2}a_{n}).
\end{eqnarray*}%
Hence, 
\begin{equation*}
\inf_{\theta \in \mathbb{R}}P_{n,\theta }\left( \theta \in C_{H,n}\right)
\leq \lim_{\theta \rightarrow \infty }p_{H,n}\left( \theta ;1,\eta
_{n},a_{n},a_{n}\right) =\Phi (n^{1/2}a_{n})-\Phi (-n^{1/2}a_{n})
\end{equation*}%
and%
\begin{equation}
\inf_{\theta \in \mathbb{R}}P_{n,\theta }\left( \theta \in E_{H,n}\right)
\leq T_{n-1}(n^{1/2}a_{n})-T_{n-1}(-n^{1/2}a_{n})\leq \Phi
(n^{1/2}a_{n})-\Phi (-n^{1/2}a_{n}),  \label{upper_bound}
\end{equation}%
the last inequality following from well-known properties of $T_{n-1}$, cf.
Lemma \ref{l_5} below. This proves the theorem in case $n^{1/2}a_{n}%
\rightarrow 0$ for $n\rightarrow \infty $.

Step 2: For every $s>0$ and $n\geq 2$ we have from (\ref{infimal_hard_asymm})%
\begin{eqnarray}
\inf_{\theta \in \mathbb{R}}P_{n,\theta }\left( \theta \in C_{H,n}\right)
&=&\inf_{\theta \in \mathbb{R}}p_{H,n}\left( \theta ;1,\eta
_{n},a_{n},a_{n}\right)  \notag \\
&=&\max \left[ \Phi (n^{1/2}a_{n})-\Phi (-n^{1/2}(a_{n}-\eta _{n})),0\right]
\label{low_bound1}
\end{eqnarray}%
and%
\begin{equation*}
\inf_{\theta \in \mathbb{R}}p_{H,n}\left( \theta ;1,s\eta
_{n},sa_{n},sa_{n}\right) =\max \left[ \Phi (n^{1/2}sa_{n})-\Phi
(n^{1/2}(-sa_{n}+s\eta _{n})),0\right] .
\end{equation*}%
Furthermore,%
\begin{eqnarray}
\inf_{\theta \in \mathbb{R}}P_{n,\theta }\left( \theta \in E_{H,n}\right)
&\geq &\int_{0}^{\infty }\inf_{\theta \in \mathbb{R}}p_{H,n}\left( \theta
;1,s\eta _{n},sa_{n},sa_{n}\right) h_{n}(s)ds  \notag \\
&=&\int_{0}^{\infty }\max \left[ \Phi (n^{1/2}sa_{n})-\Phi
(n^{1/2}(-sa_{n}+s\eta _{n})),0\right] h_{n}(s)ds  \notag \\
&=&\max \left[ \int_{0}^{\infty }\left[ \Phi (n^{1/2}sa_{n})-\Phi
(n^{1/2}(-sa_{n}+s\eta _{n}))\right] h_{n}(s)ds,0\right]  \notag \\
&=&\max \left[ T_{n-1}(n^{1/2}a_{n})-T_{n-1}(-n^{1/2}(a_{n}-\eta _{n})),0%
\right] .  \label{low_bound2}
\end{eqnarray}%
If $n^{1/2}(a_{n}-\eta _{n})\rightarrow \infty $, then the far right-hand
sides of (\ref{low_bound1}) and (\ref{low_bound2}) converge to $1$, since $%
\left\Vert \Phi -T_{n-1}\right\Vert _{\infty }\rightarrow 0$ as $%
n\rightarrow \infty $ by Polya's Theorem and since $n^{1/2}a_{n}\geq
n^{1/2}(a_{n}-\eta _{n})$. This proves the theorem in case $%
n^{1/2}(a_{n}-\eta _{n})\rightarrow \infty $.

Step 3: If $n^{1/2}\eta _{n}\rightarrow 0$, then (\ref{low_bound1}) and the
fact that $\Phi $ is globally Lipschitz shows that $\inf_{\theta \in \mathbb{%
R}}P_{n,\theta }\left( \theta \in C_{H,n}\right) $ differs from $\Phi
(n^{1/2}a_{n})-\Phi (-n^{1/2}a_{n})$ only by a term that is $o(1)$.
Similarly, (\ref{upper_bound}), (\ref{low_bound2}), the fact that $%
\left\Vert \Phi -T_{n-1}\right\Vert _{\infty }\rightarrow 0$ as $%
n\rightarrow \infty $ by Polya's theorem, and the global Lipschitz property
of $\Phi $ show that the same is true for $\inf_{\theta \in \mathbb{R}%
}P_{n,\theta }\left( \theta \in E_{H,n}\right) $, proving the theorem in
case $n^{1/2}\eta _{n}\rightarrow 0$.

Step 4: By a subsequence argument and Steps 1-3 it remains to prove the
theorem under the assumption that $n^{1/2}a_{n}$ and $n^{1/2}\eta _{n}$ are
bounded away from zero by a finite positive constant $c_{1}$, say, and that $%
n^{1/2}(a_{n}-\eta _{n})$ is bounded from above by a finite constant $c_{2}$%
, say. It then follows that $a_{n}/\eta _{n}$ is bounded by a finite
positive constant $c_{3}$, say. For given $\varepsilon >0$ set $\theta
_{n}(\varepsilon )=a_{n}(1+2c(\varepsilon )n^{-1/2})$ where $c(\varepsilon )$
is the constant given in Lemma \ref{l_3}. We then have for $s\in \lbrack
1-c(\varepsilon )n^{-1/2},1+c(\varepsilon )n^{-1/2}]$%
\begin{equation*}
sa_{n}<\theta _{n}(\varepsilon )\leq s(\eta _{n}+a_{n})
\end{equation*}%
whenever $n>n_{0}(c(\varepsilon ),c_{3})$. Without loss of generality we may
choose $n_{0}(c(\varepsilon ),c_{3})$ large enough such that also $%
1-c(\varepsilon )n^{-1/2}>0$ holds for $n>n_{0}(c(\varepsilon ),c_{3})$.
Consequently, we have (observing that $\max (0,x)$ has Lipschitz constant $1$
and $\Phi $ has Lipschitz constant $(2\pi )^{-1/2}$) for every $s\in \lbrack
1-c(\varepsilon )n^{-1/2},1+c(\varepsilon )n^{-1/2}]$ and $%
n>n_{0}(c(\varepsilon ),c_{3})$%
\begin{eqnarray*}
&&\left\vert p_{H,n}\left( \theta _{n}(\varepsilon );1,s\eta
_{n},sa_{n},sa_{n}\right) -p_{H,n}\left( \theta _{n}(\varepsilon );1,\eta
_{n},a_{n},a_{n}\right) \right\vert  \\
&=&\left\vert \max (0,\Phi (n^{1/2}sa_{n})-\Phi (n^{1/2}(-\theta
_{n}(\varepsilon )+s\eta _{n})))-\max (0,\Phi (n^{1/2}a_{n})-\Phi
(n^{1/2}(-\theta _{n}(\varepsilon )+\eta _{n})))\right\vert  \\
&\leq &\left\vert \left[ \Phi (n^{1/2}sa_{n})-\Phi (n^{1/2}(-\theta
_{n}(\varepsilon )+s\eta _{n}))\right] -\left[ \Phi (n^{1/2}a_{n})-\Phi
(n^{1/2}(-\theta _{n}(\varepsilon )+\eta _{n}))\right] \right\vert  \\
&\leq &(2\pi )^{-1/2}n^{1/2}(a_{n}+\eta _{n})\left\vert s-1\right\vert \leq
(2\pi )^{-1/2}c(\varepsilon )(a_{n}+\eta _{n})\leq (2\pi
)^{-1/2}c(\varepsilon )(c_{3}+1)\eta _{n}.
\end{eqnarray*}%
It follows that for every $n>n_{0}(c(\varepsilon ),c_{3})$%
\begin{eqnarray*}
&&\inf_{\theta \in \mathbb{R}}\int_{0}^{\infty }p_{H,n}\left( \theta
;1,s\eta _{n},sa_{n},sa_{n}\right) h_{n}(s)ds \\
&\leq &\int_{0}^{\infty }p_{H,n}\left( \theta _{n}(\varepsilon );1,s\eta
_{n},sa_{n},sa_{n}\right) h_{n}(s)ds \\
&=&\int_{1-c(\varepsilon )n^{-1/2}}^{1+c(\varepsilon )n^{-1/2}}p_{H,n}\left(
\theta _{n}(\varepsilon );1,s\eta _{n},sa_{n},sa_{n}\right) h_{n}(s)ds \\
&&+\int_{\left\{ s:\left\vert s-1\right\vert \geq c(\varepsilon
)n^{-1/2}\right\} }p_{H,n}\left( \theta _{n}(\varepsilon );1,s\eta
_{n},sa_{n},sa_{n}\right) h_{n}(s)ds \\
&=&B_{1}+B_{2}.
\end{eqnarray*}%
Clearly, $0\leq B_{2}\leq \varepsilon $ holds, cf. Lemma \ref{l_3}, and for $%
B_{1}$ we have%
\begin{eqnarray*}
&&\left\vert B_{1}-p_{H,n}\left( \theta _{n}(\varepsilon );1,\eta
_{n},a_{n},a_{n}\right) \right\vert  \\
&\leq &\left\vert \int_{1-c(\varepsilon )n^{-1/2}}^{1+c(\varepsilon
)n^{-1/2}}\left[ p_{H,n}\left( \theta _{n}(\varepsilon );1,s\eta
_{n},sa_{n},sa_{n}\right) -p_{H,n}\left( \theta _{n}(\varepsilon );1,\eta
_{n},a_{n},a_{n}\right) \right] h_{n}(s)ds\right\vert +\varepsilon  \\
&\leq &(2\pi )^{-1/2}c(\varepsilon )(c_{3}+1)\eta _{n}+\varepsilon 
\end{eqnarray*}%
for $n>n_{0}(c(\varepsilon ),c_{3})$. It follows that 
\begin{eqnarray*}
&&\inf_{\theta \in \mathbb{R}}\int_{0}^{\infty }p_{H,n}\left( \theta
;1,s\eta _{n},sa_{n},sa_{n}\right) h_{n}(s)ds \\
&\leq &p_{H,n}\left( \theta _{n}(\varepsilon );1,\eta
_{n},a_{n},a_{n}\right) +(2\pi )^{-1/2}c(\varepsilon )(c_{3}+1)\eta
_{n}+2\varepsilon 
\end{eqnarray*}%
holds for $n>n_{0}(c(\varepsilon ),c_{3})$. Now 
\begin{eqnarray*}
p_{H,n}\left( \theta _{n}(\varepsilon );1,\eta _{n},a_{n},a_{n}\right) 
&=&\max (0,\Phi (n^{1/2}a_{n})-\Phi (n^{1/2}(-\theta _{n}(\varepsilon )+\eta
_{n}))) \\
&=&\max (0,\Phi (n^{1/2}a_{n})-\Phi (n^{1/2}(-a_{n}(1+2c(\varepsilon
)n^{-1/2})+\eta _{n}))).
\end{eqnarray*}%
But this differs from $\inf_{\theta \in \mathbb{R}}P_{n,\theta }\left(
\theta \in C_{H,n}\right) =\max (0,\Phi (n^{1/2}a_{n})-\Phi
(n^{1/2}(-a_{n}+\eta _{n})))$ by at most%
\begin{eqnarray*}
&&\left\vert \Phi (n^{1/2}(-a_{n}+\eta _{n}))-\Phi
(n^{1/2}(-a_{n}(1+2c(\varepsilon )n^{-1/2})+\eta _{n}))\right\vert  \\
&\leq &(2\pi )^{-1/2}2c(\varepsilon )a_{n}\leq (2\pi )^{-1/2}2c(\varepsilon
)c_{3}\eta _{n}.
\end{eqnarray*}%
Consequently, for $n>n_{0}(c(\varepsilon ),c_{3})$ 
\begin{eqnarray*}
\inf_{\theta \in \mathbb{R}}P_{n,\theta }\left( \theta \in E_{H,n}\right) 
&=&\inf_{\theta \in \mathbb{R}}\int_{0}^{\infty }p_{H,n}\left( \theta
;1,s\eta _{n},sa_{n},sa_{n}\right) h_{n}(s)ds \\
&\leq &\max (0,\Phi (n^{1/2}a_{n})-\Phi (n^{1/2}(-a_{n}+\eta _{n})))+(2\pi
)^{-1/2}c(\varepsilon )(3c_{3}+1)\eta _{n}+2\varepsilon  \\
&=&\inf_{\theta \in \mathbb{R}}P_{n,\theta }\left( \theta \in C_{H,n}\right)
+(2\pi )^{-1/2}c(\varepsilon )(3c_{3}+1)\eta _{n}+2\varepsilon .
\end{eqnarray*}%
On the other hand, 
\begin{eqnarray*}
\inf_{\theta \in \mathbb{R}}P_{n,\theta }\left( \theta \in E_{H,n}\right) 
&=&\inf_{\theta \in \mathbb{R}}\int_{0}^{\infty }p_{H,n}\left( \theta
;1,s\eta _{n},sa_{n},sa_{n}\right) h_{n}(s)ds \\
&\geq &\int_{0}^{\infty }\inf_{\theta \in \mathbb{R}}p_{H,n}\left( \theta
;1,s\eta _{n},sa_{n},sa_{n}\right) h_{n}(s)ds \\
&=&\int_{0}^{\infty }\max (0,\Phi (n^{1/2}sa_{n})-\Phi (n^{1/2}s(-a_{n}+\eta
_{n})))h_{n}(s)ds \\
&=&\max (0,T_{n-1}(n^{1/2}a_{n})-T_{n-1}(n^{1/2}(-a_{n}+\eta _{n}))) \\
&\geq &\max (0,\Phi (n^{1/2}a_{n})-\Phi (n^{1/2}(-a_{n}+\eta
_{n})))-2\left\Vert \Phi -T_{n-1}\right\Vert _{\infty } \\
&=&\inf_{\theta \in \mathbb{R}}P_{n,\theta }\left( \theta \in C_{H,n}\right)
-2\left\Vert \Phi -T_{n-1}\right\Vert _{\infty }.
\end{eqnarray*}%
Since $\eta _{n}\rightarrow 0$ and $\left\Vert \Phi -T_{n-1}\right\Vert
_{\infty }\rightarrow 0$ for $n\rightarrow \infty $ and since $\varepsilon $
was arbitrary the proof is complete. $\ \blacksquare $

\bigskip

\noindent \textbf{Proof of Theorem \ref{adaptive_unknown}: }We prove the
result for the closed interval. Inspection of the proof together with Remark %
\ref{open}\ then gives the result for the open and half-open intervals.

Step 1: Observe that for every $s>0$ and $n\geq 2$ we have from (\ref{cov})
that%
\begin{equation*}
\lim_{\theta \rightarrow \infty }p_{A,n}\left( \theta ;1,s\eta
_{n},sa_{n},sa_{n}\right) =\Phi (n^{1/2}sa_{n})-\Phi (-n^{1/2}sa_{n}).
\end{equation*}%
Then exactly the same argument as in the proof of Theorem \ref{hard_unknown}
shows that $\inf_{\theta \in \mathbb{R}}P_{n,\theta }\left( \theta \in
C_{A,n}\right) $ as well as $\inf_{\theta \in \mathbb{R}}P_{n,\theta }\left(
\theta \in E_{A,n}\right) $ converge to zero for $n\rightarrow \infty $ if $%
n^{1/2}a_{n}\rightarrow 0$, thus proving the theorem in this case. For later
use we note that this reasoning in particular gives%
\begin{equation}
\inf_{\theta \in \mathbb{R}}P_{n,\theta }\left( \theta \in E_{A,n}\right)
\leq T_{n-1}(n^{1/2}a_{n})-T_{n-1}(-n^{1/2}a_{n})\leq \Phi
(n^{1/2}a_{n})-\Phi (-n^{1/2}a_{n}).  \label{star}
\end{equation}

Step 2: By Proposition \ref{adLASSOinf} we have for every $s>0$ and $n\geq 1$%
\begin{equation*}
\inf_{\theta \in \mathbb{R}}p_{A,n}\left( \theta ;1,s\eta
_{n},sa_{n},sa_{n}\right) =\Phi (n^{1/2}s\sqrt{a_{n}^{2}+\eta _{n}^{2}}%
)-\Phi (n^{1/2}s(-a_{n}+\eta _{n})).
\end{equation*}%
Arguing as in the proof of Theorem \ref{hard_unknown} we then have%
\begin{eqnarray}
\inf_{\theta \in \mathbb{R}}P_{n,\theta }\left( \theta \in C_{A,n}\right)
&=&\inf_{\theta \in \mathbb{R}}p_{A,n}\left( \theta ;1,\eta
_{n},a_{n},a_{n}\right)  \notag \\
&=&\Phi (n^{1/2}\sqrt{a_{n}^{2}+\eta _{n}^{2}})-\Phi (n^{1/2}(-a_{n}+\eta
_{n}))  \label{low_bound3}
\end{eqnarray}%
and%
\begin{eqnarray}
\inf_{\theta \in \mathbb{R}}P_{n,\theta }\left( \theta \in E_{A,n}\right)
&\geq &\int_{0}^{\infty }\inf_{\theta \in \mathbb{R}}p_{A,n}\left( \theta
;1,s\eta _{n},sa_{n},sa_{n}\right) h_{n}(s)ds  \notag \\
&=&T_{n-1}(n^{1/2}\sqrt{a_{n}^{2}+\eta _{n}^{2}})-T_{n-1}(n^{1/2}(-a_{n}+%
\eta _{n})).  \label{low_bound4}
\end{eqnarray}%
If $n^{1/2}(a_{n}-\eta _{n})\rightarrow \infty $, then the far right-hand
sides of (\ref{low_bound3}) and (\ref{low_bound4}) converge to $1$, since $%
\left\Vert \Phi -T_{n-1}\right\Vert _{\infty }\rightarrow 0$ as $%
n\rightarrow \infty $ by Polya's Theorem and since $n^{1/2}\sqrt{%
a_{n}^{2}+\eta _{n}^{2}}\geq n^{1/2}a_{n}\rightarrow \infty $ and $%
n^{1/2}(-a_{n}+\eta _{n})\rightarrow -\infty $. This proves the theorem in
case $n^{1/2}(a_{n}-\eta _{n})\rightarrow \infty $.

Step 3: Analogous to the corresponding step in the proof of Theorem \ref%
{hard_unknown}, using (\ref{low_bound3}), (\ref{star}), (\ref{low_bound4}),
and additionally noting that $0\leq n^{1/2}\sqrt{a_{n}^{2}+\eta _{n}^{2}}%
-n^{1/2}a_{n}\leq n^{1/2}\eta _{n}$, the theorem is proved in the case $%
n^{1/2}\eta _{n}\rightarrow 0$.

Step 4: Similar as in the proof of Theorem \ref{hard_unknown} it remains to
prove the theorem under the assumption that $n^{1/2}a_{n}\geq c_{1}>0$, $%
n^{1/2}\eta _{n}\geq c_{1}$, and that $n^{1/2}(a_{n}-\eta _{n})\leq
c_{2}<\infty $. Again, it then follows that $0\leq a_{n}/\eta _{n}\leq
c_{3}<\infty $. For given $\varepsilon >0$ set $\theta _{n}(\varepsilon
)=a_{n}(1+2c(\varepsilon )n^{-1/2})$ where $c(\varepsilon )$ is the constant
given in Lemma \ref{l_3}. We then have for $s\in \lbrack 1-c(\varepsilon
)n^{-1/2},1+c(\varepsilon )n^{-1/2}]$%
\begin{equation*}
sa_{n}<\theta _{n}(\varepsilon )
\end{equation*}%
for all $n$. Choose $n_{0}(c(\varepsilon ))$ large enough such that $%
1-c(\varepsilon )n^{-1/2}>1/2$ holds for $n>n_{0}(c(\varepsilon ))$.
Consequently, for every $s\in \lbrack 1-c(\varepsilon
)n^{-1/2},1+c(\varepsilon )n^{-1/2}]$ and $n>n_{0}(c(\varepsilon ))$ we have
from (\ref{cov}) (observing that $\Phi $ has Lipschitz constant $(2\pi
)^{-1/2}$)%
\begin{eqnarray*}
&&|p_{A,n}(\theta _{n}(\varepsilon );1,s\eta
_{n},sa_{n},sa_{n})-p_{A,n}(\theta _{n}(\varepsilon );1,\eta
_{n},a_{n},a_{n})| \\[1ex]
&\leq &(2\pi )^{-1/2}n^{1/2}\left( \left\vert s-1\right\vert
a_{n}+\left\vert \sqrt{(\theta _{n}(\varepsilon )+sa_{n})^{2}/4+s^{2}\eta
_{n}^{2}}-\sqrt{(\theta _{n}(\varepsilon )+a_{n})^{2}/4+\eta _{n}^{2}}%
\right\vert +\right. \\
&&\left. \left\vert \sqrt{(\theta _{n}(\varepsilon )-sa_{n})^{2}/4+s^{2}\eta
_{n}^{2}}-\sqrt{(\theta _{n}(\varepsilon )-a_{n})^{2}/4+\eta _{n}^{2}}%
\right\vert \right) .
\end{eqnarray*}%
We note the elementary inequality $\left\vert x^{1/2}-y^{1/2}\right\vert
\leq 2^{-1}z^{-1/2}\left\vert x-y\right\vert $ for positive $x$, $y$, $z$
satisfying $\min (x,y)\geq z$. Using this inequality with $%
z=(1-c(\varepsilon )n^{-1/2})^{2}\eta _{n}^{2}$ twice, we obtain for every $%
s\in \lbrack 1-c(\varepsilon )n^{-1/2},1+c(\varepsilon )n^{-1/2}]$ and $%
n>n_{0}(c(\varepsilon ))$%
\begin{eqnarray*}
&&|p_{A,n}(\theta _{n}(\varepsilon );1,s\eta
_{n},sa_{n},sa_{n})-p_{A,n}(\theta _{n}(\varepsilon );1,\eta
_{n},a_{n},a_{n})| \\
&\leq &(2\pi )^{-1/2}n^{1/2}|s-1|\left( a_{n}+\left[ (1-c(\varepsilon
)n^{-1/2})^{2}\eta _{n}^{2}\right] ^{-1/2}\left[ \theta _{n}(\varepsilon
)a_{n}/2+(s+1)\left( (a_{n}^{2}/4)+\eta _{n}^{2}\right) \right] \right) .
\end{eqnarray*}

Since $1-c(\varepsilon )n^{-1/2}>1/2$ for $n>n_{0}(c(\varepsilon ))$ by the
choice of $n_{0}(c(\varepsilon ))$ and since $a_{n}/\eta _{n}\leq c_{3}$ we
obtain%
\begin{eqnarray}
&&|p_{A,n}(\theta _{n}(\varepsilon );1,s\eta
_{n},sa_{n},sa_{n})-p_{A,n}(\theta _{n}(\varepsilon );1,\eta
_{n},a_{n},a_{n})|  \notag \\
&\leq &(2\pi )^{-1/2}c(\varepsilon )\left( a_{n}+2\eta _{n}^{-1}\left[
a_{n}^{2}+(5/2)((a_{n}^{2}/4)+\eta _{n}^{2})\right] \right)  \notag \\
&\leq &(2\pi )^{-1/2}c(\varepsilon )\left( c_{3}+(13/4)c_{3}^{2}+5\right)
\eta _{n}=c_{4}(\varepsilon )\eta _{n}  \label{lip_adalasso}
\end{eqnarray}%
for every $n>n_{0}(c(\varepsilon ))$ and $s\in \lbrack 1-c(\varepsilon
)n^{-1/2},1+c(\varepsilon )n^{-1/2}]$.

Now, 
\begin{align*}
& \inf_{\theta \in \mathbb{R}}\int_{0}^{\infty }p_{A,n}(\theta ;1,s\eta
_{n},sa_{n},sa_{n})h_{n}(s)ds \\
& \leq \int_{0}^{\infty }p_{A,n}(\theta _{n}(\varepsilon );1,s\eta
_{n},sa_{n},sa_{n})h_{n}(s)ds \\
& =\int_{1-c(\varepsilon )n^{-1/2}}^{1+c(\varepsilon
)n^{-1/2}}p_{A,n}(\theta _{n}(\varepsilon );1,s\eta
_{n},sa_{n},sa_{n})h_{n}(s)ds \\[1ex]
& +\int_{|s-1|\geq c(\varepsilon )n^{-1/2}}p_{A,n}(\theta _{n}(\varepsilon
);1,s\eta _{n},sa_{n},sa_{n})h_{n}(s)ds \\
& =:B_{1}+B_{2}.
\end{align*}%
Clearly, $0\leq B_{2}\leq \varepsilon $ holds by the choice of $%
c(\varepsilon )$, see Lemma \ref{l_3}. For $B_{1}$ we have using (\ref%
{lip_adalasso})%
\begin{eqnarray*}
&&|B_{1}-p_{A,n}(\theta _{n}(\varepsilon );1,\eta _{n},a_{n},a_{n})| \\
&\leq &\int_{1-c(\varepsilon )n^{-1/2}}^{1+c(\varepsilon
)n^{-1/2}}|p_{A,n}(\theta _{n}(\varepsilon );1,s\eta
_{n},sa_{n},sa_{n})-p_{A,n}(\theta _{n}(\varepsilon );1,\eta
_{n},a_{n},a_{n})|h_{n}(s)ds+\varepsilon  \\
&\leq &c_{4}(\varepsilon )\eta _{n}+\varepsilon 
\end{eqnarray*}%
for $n>n_{0}(c(\varepsilon ))$. It follows that 
\begin{eqnarray*}
&&\inf_{\theta \in \mathbb{R}}\int_{0}^{\infty }p_{A,n}(\theta ;1,s\eta
_{n},sa_{n},sa_{n})h_{n}(s)ds \\
&\leq &p_{A,n}(\theta _{n}(\varepsilon );1,\eta
_{n},a_{n},a_{n})+c_{4}(\varepsilon )\eta _{n}+2\varepsilon 
\end{eqnarray*}%
holds for $n>n_{0}(c(\varepsilon ))$. Furthermore, the absolute difference
between $p_{A,n}(\theta _{n}(\varepsilon );1,\eta _{n},a_{n},a_{n})$ and $%
\inf_{\theta \in \mathbb{R}}P_{n,\theta }\left( \theta \in C_{A,n}\right) $
can be bounded as follows: Using Proposition \ref{adLASSOinf}, (\ref{cov}),
observing that $\Phi $ has Lipschitz constant $(2\pi )^{-1/2}$, and using
the elementary inequality noted earlier twice with $z=\eta _{n}^{2}$ we
obtain%
\begin{eqnarray*}
&&\left\vert p_{A,n}(\theta _{n}(\varepsilon );1,\eta _{n},a_{n},a_{n})-\Phi
\left( n^{1/2}\sqrt{a_{n}^{2}+\eta _{n}^{2}}\right) +\Phi \left(
n^{1/2}(-a_{n}+\eta _{n})\right) \right\vert  \\
&\leq &(2\pi )^{-1/2}n^{1/2}\left\vert -a_{n}c(\varepsilon )n^{-1/2}+\sqrt{%
a_{n}^{2}(1+c(\varepsilon )n^{-1/2})^{2}+\eta _{n}^{2}}-\sqrt{a_{n}^{2}+\eta
_{n}^{2}}\right\vert  \\
&&+(2\pi )^{-1/2}n^{1/2}\left\vert \sqrt{(a_{n}c(\varepsilon
)n^{-1/2})^{2}+\eta _{n}^{2}}-\sqrt{(a_{n}c(\varepsilon )n^{-1/2}+\eta
_{n})^{2}}\right\vert  \\
&\leq &(2\pi )^{-1/2}\left( 2a_{n}c(\varepsilon )+(2\eta
_{n})^{-1}a_{n}^{2}(2c(\varepsilon )+c(\varepsilon )^{2}n^{-1/2})\right)  \\
&\leq &(2\pi )^{-1/2}\left( 2c_{3}c(\varepsilon
)+2^{-1}c_{3}^{2}(2c(\varepsilon )+c(\varepsilon )^{2})\right) \eta
_{n}=c_{5}(\varepsilon )\eta _{n}.
\end{eqnarray*}%
Consequently, for $n>n_{0}(c(\varepsilon ))$ 
\begin{align*}
& \inf_{\theta \in \mathbb{R}}\int_{0}^{\infty }p_{A,n}(\theta ;1,s\eta
_{n},sa_{n},sa_{n})h_{n}(s)ds \\
& \leq \Phi (n^{1/2}\sqrt{a_{n}^{2}+\eta _{n}^{2}})-\Phi
(n^{1/2}(-a_{n}+\eta _{n})) \\
& +\left( c_{4}(\varepsilon )+c_{5}(\varepsilon )\right) \eta
_{n}+2\varepsilon .
\end{align*}%
On the other hand, 
\begin{align*}
& \inf_{\theta \in \mathbb{R}}\int_{0}^{\infty }p_{A,n}(\theta ;1,s\eta
_{n},sa_{n},sa_{n})h_{n}(s)ds \\
& \geq \int_{0}^{\infty }\inf_{\theta \in \mathbb{R}}p_{A,n}(\theta ;1,s\eta
_{n},sa_{n},sa_{n})h_{n}(s)ds \\
& =\int_{0}^{\infty }\left[ \Phi (n^{1/2}s\sqrt{a_{n}^{2}+\eta _{n}^{2}}%
)-\Phi (n^{1/2}s(-a_{n}+\eta _{n}))\right] h_{n}(s)ds \\
& =T_{n-1}(n^{1/2}\sqrt{a_{n}^{2}+\eta _{n}^{2}})-T_{n-1}(n^{1/2}(-a_{n}+%
\eta _{n})) \\
& \geq \Phi (n^{1/2}\sqrt{a_{n}^{2}+\eta _{n}^{2}})-\Phi
(n^{1/2}(-a_{n}+\eta _{n}))-2\Vert \Phi -T_{n-1}\Vert _{\infty }.
\end{align*}%
Since $\eta _{n}\rightarrow 0$ and $\left\Vert \Phi -T_{n-1}\right\Vert
_{\infty }\rightarrow 0$ for $n\rightarrow \infty $ and since $\varepsilon $
was arbitrary the proof is complete. $\ \blacksquare $

\bigskip

\begin{lemma}
\label{l_3}Suppose $\sigma =1$. Then for every $\varepsilon >0$ there exists
a $c=c(\varepsilon )>0$ such that 
\begin{equation*}
\int_{\max (0,1-cn^{-1/2})}^{1+cn^{-1/2}}h_{n}(s)ds\geq 1-\varepsilon
\end{equation*}%
holds for every $n\geq 2$.
\end{lemma}

\begin{proof}
By the central limit theorem and the delta-method we have that $n^{1/2}(\hat{%
\sigma}-1)$ converges to a normal distribution. It follows that $n^{1/2}(%
\hat{\sigma}-1)$ is (uniformly) tight. In other words, for every $%
\varepsilon >0$ we can find a real number $c>0$ such that for all $n\geq 2$
holds%
\begin{equation*}
\Pr \left( \left\vert n^{1/2}(\hat{\sigma}-1)\right\vert \leq c\right) \geq
1-\varepsilon .
\end{equation*}
\end{proof}

\begin{lemma}
\label{l_5} Suppose $n\geq 2$ and $x\geq y\geq 0$. Then%
\begin{equation*}
T_{n-1}(x)\leq \Phi (x)
\end{equation*}%
and%
\begin{equation*}
T_{n-1}(x-y)-T_{n-1}(-x-y)\leq \Phi (x-y)-\Phi (-x-y).
\end{equation*}
\end{lemma}

\begin{proof}
The first claim is well-known, see, e.g., Kagan and Nagaev (2008). The
second claim follows immediately from the first claim, since by symmetry of $%
\Phi $ and $T_{n-1}$ we have 
\begin{eqnarray*}
&&\Phi (x-y)-\Phi (-x-y)-\left( T_{n-1}(x-y)-T_{n-1}(-x-y)\right) \\
&=&\left[ \Phi (x-y)-T_{n-1}(x-y)\right] +\left[ \Phi (x+y)-T_{n-1}(x+y)%
\right] \geq 0.
\end{eqnarray*}
\end{proof}

\end{document}